\documentclass[10pt]{article}
\usepackage[T1]{fontenc}
\usepackage[utf8]{inputenc}
\usepackage{amsmath,amssymb,latexsym,soul,cite,mathrsfs}
\usepackage{amsfonts}
\usepackage{amssymb}
\usepackage{graphicx}
\usepackage{verbatim}
\usepackage{mathrsfs}
\usepackage{upref,amsthm,amsxtra,exscale}
\usepackage{soul}
\usepackage{nameref}
\usepackage{varioref}
\usepackage[colorlinks=true,urlcolor=blue,
citecolor=red,linkcolor=blue,linktocpage,pdfpagelabels,
bookmarksnumbered,bookmarksopen,pagebackref]{hyperref}
\usepackage{cleveref}
\usepackage{cite}
\usepackage{fullpage}
\usepackage[dvipsnames]{xcolor}
\usepackage{esint}
\newtheorem{theorem}{Theorem}[section]
\newtheorem{corollary}[theorem]{Corollary}
\newtheorem{remark}[theorem]{Remark}

\newtheorem{lemma}[theorem]{Lemma}
\newtheorem{proposition}[theorem]{Proposition}
\newtheorem{definition}[theorem]{Definition}

\numberwithin{equation}{section}

\def\r{\mathbb{R}}

\def\s1{\mathbb{S}^1}

\def\eps{\varepsilon}

\def\cH^G{\mathcal{H}}

\def\bar{\overline}

\def\i{\mathrm{i}}

\author{Romildo Lima\footnote{R. Lima was partially supported by CNPq/Brazil 306.411/2022-9 and 176.596/2023-2.}\\
Unidade Acadêmica de Matemática, Universidade Federal de Campina Grande\\ 58429-970, Campina Grande - PB, Brazil\\	
Mayra Soares,\footnote{M. Soares was partially supported by CNPq/Brazil 303.154/2025-0 and 152.786/2025-2.} \\
Departamento de Matemática, Universidade de Brasília \\70910, Brasília-DF, Brazil\\
 Leandro Tavares \\(corresponding author/e-mail: leandrolstav@gmail.com)\\
Centro de Matemática, Computação e Cognição, Universidade Federal do ABC\\ 63048-080, Santo André - SP, Brazil}

\title{Regularity Estimates for a Nonhomogeneous Free Boundary Interface Problem}

\begin{document}
		\maketitle
	
\begin{abstract}
In this work, we establish regularity properties for minimizers of a novel nonhomogeneous interface problem with a free boundary, formulated in Orlicz–Sobolev spaces. The nonlinear terms fulfill Lieberman-type structural assumptions with subcritical growth. The principal difficulties stem from the nonhomogeneous framework and the lack of suitable representation formulas. Our findings extend recent contributions in the field and provide new insights even in the particular case of the Laplacian.


\end{abstract}

%\tableofcontents
\section{Introduction and Main Results}\label{intro}
\hfill

In this paper, which was motivated by \cite{LQT,zhengt,ZFZ,ZZZ,bianc,piment}, we are interested in the  minimization problem 
\begin{equation}\label{func}
 J(v):= \int_{\Omega} \Phi(|\nabla v|) \ \text{d}x + \int_{\Omega_1} W(v)  + gv \ \text{d}x +\int_{\Omega_2} \overline{W}(v) +\overline{g}v \ \text{d}x - \int_{\Gamma}  hv - Z(v) \ \text{d}H^{N-1} \rightarrow {\textrm{min}}, 
\end{equation}
 which arises in the study of transmission problems, where $|\cdot|$ denotes the usual norm in $\mathbb{R}^{N}, \ N\geq 3$, $\Phi(t) = \int_{0}^{t}\phi(s)\text{d}s$ and
 \begin{equation}\label{parametersss}
 \begin{aligned}
 W(v)&:=&\theta' F(v^{+}) +\widetilde{\theta} F(v^{-}) + \lambda^{\prime}\chi_{\{v>0\}} +\widetilde{\lambda}\chi_{\{v\leq 0\}}
 ,\\
 \overline{W}(v) &:= &  \sigma'\overline{F}(v^{+}) +\widetilde{\sigma} \overline{F}(v^{-}) + \beta^{\prime}\chi_{\{v>0\}} +\widetilde{\beta}\chi_{\{v\leq0\}}
 ,\\
  Z(v)&:=&  \rho'Y(v^{+}) + \widetilde{\rho}Y(v^{-}) + \gamma^{\prime}\chi_{\{v>0\}} +\widetilde \gamma\chi_{\{v\leq0\}},	
 \end{aligned}
\end{equation}
with $v^{\pm}:= \max\{ \pm v,0\}$, $\theta',\widetilde{\theta},\sigma',\widetilde{\sigma},\rho',\widetilde{\rho},\widetilde{\lambda},\widetilde{\beta}, \widetilde\gamma, \lambda^{\prime}, \beta^{\prime}, \gamma^{\prime} \geq 0$  are  constants satisfying $\widetilde \lambda\leq \lambda'$, $\widetilde \beta\leq \beta'$, $\widetilde \gamma\leq \gamma'$, and $\Phi^{\prime}(t) = \phi(t),$ for all $t>0,$
over the set 
\begin{equation*}
W^{1,\Phi}_{\varphi} (\Omega):=  \left\{v \in W^{1,\Phi}(\Omega) ; v-\varphi \in W^{1,\Phi}_{0}(\Omega)\right\},
\end{equation*}
where $W^{1,\Phi}(\Omega)$ denotes the Orlicz Space defined by $\Phi,$ $\varphi \in W^{1,\Phi}(\Omega)\cap L^{\infty}(\Omega),$  and $H^{N-1}$ is  the  Hausdorff measure in the surface integral. It is a fact that we are considering $\Omega_1 \subset \subset\Omega \subset \mathbb{R}^N$  bounded domains with smooth boundaries and $\Gamma:= \partial \Omega_1$ and $\Omega_2 := \Omega\setminus  \overline{\Omega_1}$. In addition, $\phi,F,\overline{F}, Y \in C([0,\infty);[0,\infty))\cap C^1((0,\infty);(0,\infty))$ with $\phi(0)=F(0)=\overline{F}(0)=Y(0)=0$ satisfying, for positive constants $\delta_0,g_{0}\in\r$
\begin{align}\label{condphi}
	0<{\delta_0}\leq \frac{t\phi'(t)}{\phi(t)}\leq g_{0},\ \text{for all} \   t>0,
\end{align}
with $1+\delta_0 <N,$ and
\begin{equation}\label{hipF}
1+\theta_0\leq \frac{F^{\prime}(t)t}{F(t)}\leq 1+f_{0}, 	1+\overline{\theta}_0\leq \frac{{(\overline{F})}^{\prime}(t)t}{\overline{F}(t)}\leq 1+\overline{f}_{0}, 	1+\widetilde{\theta}_0\leq \frac{Y^{\prime}(t)t}{Y(t)}\leq 1+\widetilde{f}_{0}, \ \text{for all} \   t >0,
\end{equation}
where $\theta_0,f_0,\overline{\theta}_0,\overline{f}_{0},\widetilde{\theta}_0,\widetilde{f}_{0}\in\r$, $F^{\prime}(t)=f(t), {\overline{F}}^{\prime}(t)=\overline{f}(t), Y^{\prime}(t)=y(t)$, for all $t>0$,  and
\begin{equation*}
0<1+\theta_0\leq 1+f_0 <(1+\delta_0)^{*},
\quad 
0<1+\overline{\theta}_0\leq 1+\overline{f}_0 <(1+\delta_0)^{*}, \quad
0<1+\widetilde{\theta}_0\leq 1+\widetilde{f}_0 <(1+\delta_0)_{*},
\end{equation*}
considering that $(1+\delta_0)_{*}:=\frac{(N-1)(1+\delta_0)}{N - (1+\delta_0)}$ and $(1+\delta_0)^* := \frac{N(1+\delta_0)}{N-(1+\delta_0)}$ is the Sobolev exponent of $1+\delta_0$. 
It is important to note that
\begin{align}\label{generalhyp}
	g \in L^{q}(\Omega_1),\overline{g} \in L^{\overline{q}}(\Omega_2) \ \text{with} \  q,\overline{q} \geq N, \ \text{and} \ h \in L^{\frac{m}{m-1}}(\Gamma),   \text{ for some} \ 1\leq  m < (1+\delta_0)_{*}.
\end{align}

The energy functional \eqref{func} may be seen as a generalization of a class of functionals that have been extensively studied in the literature since the 1980's, a period characterized by significant advances in the regularity theory of minimizers for variational problems involving free boundaries and their associated geometric structures. For example, when ($\lambda^{\prime}, \beta^{\prime} > 0$) and all the other parameters in \eqref{parametersss} are zero, our  problem refers to the remarkable work  \cite{altcaf}, by  Alt and Caffarelli, who examined the one-phase free boundary problem -- closely connected to the classical obstacle problem and Bernoulli-type problems -- and proved the Lipschitz regularity of minimizers of the corresponding functional
\begin{equation*}
	\mathcal{S}(v):=\int_{\Omega}  |\nabla v|^2 + Q(x) \chi_{\{v>0\}}  {\textrm{d}}x,\mbox{ for } v \in  K_{\varphi} := \left\{ v \in H^1(\Omega) ; v - \varphi \in H_0^1(\Omega) \right\}.
\end{equation*}
They considered
$0 \leq \varphi \in H^1(\Omega) \cap L^\infty(\partial\Omega)$ and $0 < \lambda \leq Q(x) \leq \Lambda$, which  provided the foundation for subsequent researches in this field. For instance, in \cite{acf}, Alt, Caffarelli, and Friedman analyzed the two-phase cavity problem, further enhancing the understanding of free boundary problems.
Similarly, we highlight the contribution of Alt and Phillips in  \cite{altp}, which is related to the case $\theta^{\prime},\sigma^{\prime} > 0$, with the other parameters vanishing in \eqref{parametersss}, and  $F(t),\overline{F}(t) = t^{\gamma}, t \geq 0$, where $0 < \gamma < 2$. Their work focuses on studying the minimizers of a broad class of functionals that includes the case  
\[
\mathcal{S}_{\gamma}(v) := \int_{\Omega} \frac{|\nabla v|^2}{2} + (v^{+})^{\gamma} \, \mathrm{d}x,
\]  
over an appropriate subset of $W^{1,2}(\Omega)$. These minimizers are used to model the distribution of a gas whose density corresponds to the minimizer within a porous catalyst pellet $\Omega$.

Subsequently, these previous references inspired several generalizations to the nonhomogeneous setting. For example, Danielli and Petrosyan \cite{dp} extended the work of \cite{altcaf} to the $p-$Laplacian case. Additionally, Karakhanyan \cite{karak} studied the regularity of minimizers within a suitable subspace of $W^{1,p}(\Omega)$ for the functional given by
$$
\mathcal{E}_p(v, \Omega) = \int_{\Omega} \left[ |\nabla v|^p + \lambda_1^p \, \chi_{\{v > 0\}} + \lambda_2^p \, \chi_{\{v \leq 0\}} \right], \quad 1 < p < \infty,
$$
where $\lambda_1, \lambda_2$ are positive constants such that $\Lambda = \lambda_2^p - \lambda_1^p < 0$, which  is related to the case where $\lambda', \widetilde{\lambda} > 0$ and all other parameters vanish in \eqref{func}, and that arise in models involving equilibrium configurations for heat or electrostatic energy optimization.

 In the case of one phase problems involving  operators that generalize the Ladyzhenskaya–Uraltseva operators, the first reference is Martínez and Wolanski \cite{MW}, which considers an optimization problem that generalizes \cite{altcaf} in an Orlicz-Sobolev setting. It was proved that every minimizer  is locally Lipschitz continuous and that the free boundary is a regular surface. We also highlight the works of Challal and Lyaghfouri \cite{CL1,CL2}, which considered the obstacle problem case. They established the porosity of the free boundary in \cite{CL1} and proved that the free boundary is locally of finite Hausdorff measure in the reference \cite{CL2}.
 
We also highlight \cite{zhengt}, which complements \cite{MW} and is related to the case where $\lambda^{\prime},\beta^{\prime},\theta^{\prime},\sigma^{\prime}>0$, while the remaining parameters in \eqref{parametersss} are set to zero. In that work, with certain measurable functions $q$, $h$, and a nonnegative constant $\lambda$, the authors established regularity results for the free boundary problem  
$$
\mathcal{J}(v) := \int_{\Omega}  \Phi(|\nabla v|) + qF(v^+) + h v^+ +  \lambda \chi_{\{v>0\}}  {\textrm{d}}x \to \min,
$$
where the functional is defined in the class of functions $v \in  W^{1,\Phi}(\Omega) $ with $ v - \varphi \in W^{1,\Phi}_0(\Omega) $ for a fixed bounded function $ \varphi \in W^{1,\Phi}(\Omega),$
where  $F$ satisfy the structural conditions of Liebermann’s type given in \eqref{hipF}, and exhibit subcritical  growth with $h $ belonging to a suitable Lebesgue space. The local $C^{1,\alpha}$ and Log-Lipschitz regularity for minimizers were obtained in the subcritical case for $\lambda = 0$ and $\lambda \geq 0$, respectively. The authors established growth rates near the free boundary for non-negative minimizers in the subcritical case for $\lambda \geq 0$, as well as the optimal growth rate of non-negative minimizers in the $p$-Laplacian case. Moreover, they proved the local Lipschitz continuity of non-negative minimizers for $\lambda > 0$, assuming the growth condition $F(t) \leq 1 + \Phi (t)$ for $t \geq 0$.

Regarding free boundary problems involving two phases for the $p$-Laplacian operator, we refer to the work of Leit\~{a}o, de Queiroz, and Teixeira \cite{LQT}. Their analysis dealt with  the case where $
		\theta', \sigma', \lambda', \beta', \widetilde{\theta}, \widetilde{\sigma}, \widetilde{\lambda}, \widetilde{\beta} > 0,$
	while all other parameters in \eqref{parametersss} vanish. The authors established sharp regularity results for a class of problems including both Prandtl--Batchelor-type and obstacle-type problems. More precisely, they proved regularity estimates for the associated minimization problem
\begin{equation*}
	J_{\gamma} (v):=	\int_{\Omega}  |\nabla v|^p + F_{\gamma}(v) + f  v \  {\textrm{d}}x \to \min,
\end{equation*}
for  $ v \in W^{1,p}(\Omega) + \varphi , \varphi \in W^{1,p}(\Omega)\cap L^{\infty}(\partial \Omega)$,$ \ \gamma \in [0,1] $, and 
\begin{equation*}
	F_{\gamma}(v) := \lambda_+ (v^+)^{\gamma} + \lambda_- (v^-)^{\gamma}, \gamma \in (0,1], \quad \mbox{ where }  0 \leq \lambda_- < \lambda_+ ,
\end{equation*}
\begin{equation*}
	F_0(v) := \lambda_+ \chi_{\{v>0\}} + \lambda_- \chi_{\{v\leq 0\}},
\end{equation*}
and  it was  investigated the asymptotic behavior as \( \gamma \to 0 \).

Motivated by \cite{LQT}, Zheng, Zhang, and Zhao in \cite{ZZZ} studied the minimizers of the functional  
$$
\mathcal{J}_\gamma (v) := \int_{\Omega} \Phi(|\nabla v|) + \lambda_+ (v^+)^\gamma + \lambda_- (v^-)^\gamma + f v \  \textrm{d}x,
$$ 
defined in the class of functions $ W^{1,\Phi}(\Omega) $ with $ u - \varphi \in W^{1,\Phi}_0(\Omega) $ for a fixed bounded function $ \varphi \in W^{1,\Phi}(\Omega) $, and $f$ is a measurable function satisfying certain conditions. Several regularity results were obtained for the minimizers of $ \mathcal{J}_{\gamma} $. In the case \( 0 < \gamma \leq 1 \), it was shown that every minimizer is \( C^{1,\alpha}_{\rm loc} \)-continuous. For \( \gamma = 0 \), the authors established local log-Lipschitz continuity for any minimizer. As in \cite{LQT}, the authors  considered  the asymptotic behavior as \( \gamma \to 0 \) and derived sharp geometric estimates near the free boundary in the case \( \gamma = 0 \). Improved results for the case $0 < \gamma \leq 1 $ were obtained by Zheng, Feng, and Zhao \cite{ZFZ}.

Transmission problems, whose mathematical formulation involves partial differential equations in a domain split into sub-regions, arise in several areas of physics and engineering. As pointed out in \cite{borsuk},  such problems arise  in  electrodynamics of solid media in the consideration of  electromagnetic processes in ferromagnetic materials with varying dielectric constants, and  in the study of  composite materials, vibrating folded membranes, composite plates, folded plates,  junctions in elastic multi-structures.

A notable aspect of transmission problems is the fact that  the partial differential equations governing each sub-region of the domain may differ, influencing the behavior of solutions at the fixed interface. Thus, the assumptions concerning the structure of operators in the regions and at the interfaces play a pivotal  role in the analysis of the problem. The mathematical interest in transmission problems dates back to the 1930s, notably with the  works of Picone \cite{Picone2,Picone}, which was motivated by elasticity theory, and Lions \cite{lions} that  considered the  existence and uniqueness of solutions. Schechter \cite{schechter} later extended the study to high-order elliptic equations with general boundary conditions. Further contributions can be found in  Il’in \cite{iin} and  \v{S}i\v{s}marev \cite{iin2}, which  applied potential methods to establish classical solvability for Dirichlet and Neumann problems with  general linear elliptic second-order operators, and discontinuous coefficients in nondivergence form. Other significant references include  \cite{oleinik,camp,camp2,camp3,iin,LU,sheftel,stampacchia,borsuk2}. For a detailed  overview of the historical aspects of transmission problems, we cite  the recent reference \cite{borsuk}.

Regarding regularity results for interface problems, one of the earliest significant contributions to the H\"{o}lder continuity of the gradient of solutions is the classical paper due to  Li and Vogelius \cite{livo}. In such reference  it was obtained  $W^{1,\infty}$ and piece-wise $C^{1,\alpha}$ estimates for solutions to divergence-form elliptic equations with piece-wise H\"{o}lder continuous coefficients,   whose  domain is split into a finite number of pairwise disjoint subdomains, with the operator potentially exhibiting a ``blow-up''  near the boundary of these subdomains.

Recently, Caffarelli, Soria-Carro and Stinga \cite{caf}, which focuses in the minimal requirements of the  geometry of the interface, considered  the transmission problem which consists in finding a continuous function $v: \overline{\Omega}\rightarrow \mathbb{R}$  of the  transmission problem
\begin{equation}\label{cafp}
\begin{cases}
	\Delta v_1 = 0 \quad \text{in } \Omega_1, \quad v = 0 \quad \text{on } \partial\Omega, \\
	\Delta v_2 = 0 \quad \text{in } \Omega_2, \quad	v_1 = v_2 \quad \text{on } \Gamma, \\
	\dfrac{\partial v_1 }{\partial \nu  } - \dfrac{\partial v_2 }{\partial \nu  }= h \quad \text{on } \Gamma,
\end{cases}
\end{equation}
with   $\Omega$, $\Omega_1$, $\Omega_2$ and $\Gamma$ as previously described, but only with the requirement that $\Gamma $ is a  $C^{1,\zeta}$ manifold,  $\nu$ denotes  the unit normal vector on $ \Gamma $ that is interior to $\Omega_1 $, $h\in C^{0,\zeta}(\Gamma)$ is a non-negative function  $ v_1 = v|_{\Omega_1}$ and $ v_2 = v|_{\Omega_2} $. The authors proved \eqref{cafp} admits a unique distributional solution in $C(\overline{\Omega}),$ and  that such solution is Log-Lipschitz in $\overline{\Omega}.$ Moreover, it is also proved that \eqref{cafp} admits a unique classical solution, and the authors showed that it lies in $ C^{1,\zeta} $ from each side up to \( \Gamma \). Such notable results  were obtained by means  of a representation formula  for harmonic functions, the maximum principle, and  new  geometric stability type arguments. An alternative proof of the $ C^{1,\zeta} $ estimate  of \cite{caf}, was obtained in \cite{dong}. Such a proof does not require the application of representation formulas for harmonic functions or maximum principles. Moreover, it applies to a broader class of elliptic systems with variable coefficients.

In the nonhomogeneous setting, where in general representation formulas are not available, an interface problem was recently considered in \cite{piment,bianc}. Specifically, in such references  it was  obtained the existence of    minimizers for the energy functional 
\begin{equation}\label{pimentf}
 I(v):= \int_{\Omega} \Phi(|\nabla v|)\text{d}x - \int_{\Gamma} h v \ \text{d}H^{N-1}, v \in W^{1,\Phi}_{0}(\Omega),
\end{equation}
where $h$ is a measurable function satisfying certain conditions, under the structural condition
	\begin{equation}\label{ed}
	\left\langle \frac{\phi(\xi)}{|\xi|} \xi - \frac{\phi(\eta)}{|\eta|} \eta ,  \xi-\eta\right\rangle \geq C |\xi-\eta|^{p}, \mbox{ for } p>2, \  \delta_0 \geq 1,
\end{equation}
 and some  constant $C>0$,  for  all $\xi,\eta \in \mathbb{R}^N  - \{0\},$  where $\langle \cdot, \cdot \rangle$ denotes the usual inner product in $\mathbb{R}^N$. It was showed that the minimizers coincides with the weak solutions of the interface problem given by

\begin{equation*}
	\begin{cases}
		\Delta_{\Phi}v_1  = 0 \quad \text{in } \Omega_1, \\
		\Delta_{\Phi} v_2  = 0 \quad \text{in } \Omega_2, \\
		v = 0 \quad \text{on } \partial\Omega, \\
		\dfrac{\partial v_1 }{\partial \nu  } - \dfrac{\partial v_2 }{\partial \nu  }= h \quad \text{on } \Gamma,
	\end{cases}
\end{equation*}
with $\Omega$, $\Omega_1$, $\Omega_2$, $\Gamma$, $\Phi$, $\phi$, $\nu$, $\Gamma \in C^{1}$, $v_1$, $v_2$ as before, and $\Delta_{\Phi} v_i := \textrm{div}\left( \tfrac{\phi(|\nabla v_i|)}{|\nabla v_i|} \nabla v_i\right),i=1,2. $ In their first result, motivated by Serrin \cite{serrin}, the authors established local boundedness and local estimates for the $L^p-$norm of the gradient of a minimizer of \eqref{pimentf}, assuming $ \phi(t) = t^{p-1} $ and $ h \in L^{\infty}(\Gamma) $. When $ h \in L^{\infty}(\Gamma) $, they further showed that minimizers exhibit local BMO regularity by  comparing them with $\Phi$-harmonic functions. Additionally, Log-Lipschitz estimates were obtained across the transmission interface.  In the case of unbounded interface data, under the assumption that \( h \in W^{1, p' + \varepsilon}(\Omega) \) for some suitable \( \varepsilon > 0 \), the authors showed, among other results, that minimizers belong to a certain Campanato space.

Therefore, based on the previous comments we propose the study of minimizers of the functional \eqref{func}. The purpose of this manuscript is to complement the results obtained in the recent reference \cite{piment}. Below we describe precisely the results that were obtained, whose proofs are based on the ideas of \cite{bianc,LQT,ZZZ,piment,queirozt,tavares,ZFZ,zhengt}.
\subsection{Presentation of the Results and Examples}
\hfill

The first result obtained consists in providing the existence of minimizers for the functional \eqref{func}. { This theorem covers a broad family of nonlinearities and a general transmission term, thus extending and unifying several earlier results in the literature.} 
\begin{theorem} \label{existence-boundness}
	{Suppose that \eqref{condphi}, \eqref{hipF}, and \eqref{generalhyp} hold. Then, the functional ${J}$ admits at least one minimizer.}
\end{theorem}
{In particular, Theorem \ref{existence-boundness} extends the existence result obtained in \cite{piment} by considering \(\lambda^{\prime} = \beta^{\prime} = \gamma^{\prime}= \widetilde{\lambda}=\widetilde{\beta}=\widetilde\gamma=0\) in \eqref{func}. Moreover, the proof presented here employs a different method compared to the reference  \cite{piment}. 
An important aspect of the proof of Theorem \ref{existence-boundness}, which is  based on the ideas developed in \cite{LQT,ZZZ,ZFZ,zhengt}, is the presence of general transmission terms and the  other nonlinear terms. Such terms introduce additional mathematical difficulties compared to \cite{LQT,piment, ZZZ,ZFZ,zhengt}, requiring a refinement of some estimates of such references to handle the complexities that arise in the considered setting.}

{We now present results in the case of a bounded interface data, that is, $h \in L^{\infty}(\Gamma).$ The result below assures the $L^{\infty}$ regularity of the minimizers, and uniform estimates. Unlike \cite[Theorem 1]{piment}, Theorem \ref{unifbound}, provides uniform estimates in $L^{\infty}(\Omega),$ and $W^{1,\Phi}(\Omega)$ for the minimizers of \eqref{func} in the entire domain, thereby complementing \cite[Theorem 1]{piment} that considered local estimates in the mentioned spaces.}

\begin{theorem}\label{unifbound} Suppose that \eqref{condphi}, \eqref{hipF},  \eqref{generalhyp} hold and that $h \in L^{\infty}(\Gamma).$ Then $v \in L^{\infty}(\Omega).$   Moreover, there is a universal constant $K_\varphi$,  such that
	\begin{equation}\label{lpart}
		\| v\|_{L^{\infty}(\Omega)}+\|v\|_{W^{1,\Phi}(\Omega)} \leq K_\varphi,
	\end{equation}
	holds true for all minimizer $v$ of the functional $J$.
\end{theorem}

{ The proof of such a result is based on a nontrivial adaptation of the ideas in \cite{LQT,ZZZ,ZFZ,zhengt} and employs a different approach from the one used in the reference  \cite[Theorem 1]{piment}. In particular, the presence of the general nonlinearities and transmission terms  needs different estimates compared to those in \cite{LQT, piment, ZZZ, ZFZ, zhengt}, in order to derive the uniform estimates throughout the domain, including across the transmission interface, given in Theorem \ref{unifbound}.}

Before presenting the next result, we introduce a notation. 
Consider $\Omega \subset \mathbb{R}^N$  an open domain and  $\zeta \geq 0,$ and the Morrey-type space
$$ M^{\Phi , \zeta}(\Omega):= \left\{v \in L^{\Phi}(\Omega) ;  \displaystyle\sup_{x\in \Omega, R>0} R^{-\zeta} \int_{\Omega \cap {B}_R(x)} \Phi(|v|) \text{d}x < \infty \right\},  $$
where    $L^{\Phi}(\Omega)$ denotes the usual Orlicz space ( see Section \ref{suborlicz} for more details on such spaces), and $\Phi$ satisfies \eqref{condphi}. It will be denoted by $M^{\Phi , \zeta}_{\rm{loc} }(\Omega)$ the set of functions  that belongs to $M^{\Phi , \zeta}(\Omega^{\prime})$ for all open set $\Omega^{\prime}\subset \subset \Omega.$ For more details regarding Morrey spaces, see for example \cite{giu}.
In the following result,  Morrey-type estimates for the minimizers of the functional \eqref{func} were obtained.

{
	In the next theorem, where Theorem \ref{unifbound} played a pivotal role, integrability results were obtained for the gradient of the minimizers of the problem  \eqref{func}, an issue that was not considered in \cite{piment}. Specifically, Morrey-type estimates were obtained in Theorem \ref{morreyg} for these functions.}

\begin{theorem}\label{morreyg}  Suppose that \eqref{condphi} , \eqref{hipF}, \eqref{generalhyp}  are in force, and also that  $h \in L^{\infty}(\Gamma).$  The assertions below hold for a minimizer $v$ of the functional \eqref{func}. If $y  = Y^{\prime}
	$ is bounded or nondecreasing, and
\begin{itemize}
		\item[(a)]  $\gamma^{\prime}>0$ or $\widetilde\gamma>0,$  then	 $|\nabla v| \in M^{\Phi , \zeta}_{\rm{loc}}(\Omega)$ for all $0\leq \zeta \leq  N-1;$
		\item[(b)]  $\gamma^{\prime} = \widetilde\gamma=0,$ then	 $|\nabla v| \in M^{\Phi , \zeta}_{\rm{loc}}(\Omega)$ for all $0\leq \zeta < N.$
\end{itemize}		
	On the other hand, 	if $y = Y^{\prime} $ is nonincreasing with $-1< \widetilde{\theta}_0 \leq \widetilde{f}_0 <0$ or $-1 < \theta_0 < 0 \leq \widetilde{f}_0$, and
	\begin{itemize}
	 \item[(c)]   $\gamma^{\prime}>0 $ or $\widetilde\gamma >0,$ then	$|\nabla v| \in M^{\Phi , \zeta}_{\rm{loc}}(\Omega)$ for all $0\leq  \zeta\leq N-1;$
		\item[(d)]  $\gamma^{\prime} = \widetilde\gamma=0,$ then $|\nabla v| \in M^{\Phi , \zeta}_{\rm{loc}}(\Omega)$ for all $0\leq  \zeta <\, N$.
		\end{itemize}
		Finally, if $y = Y^{\prime} $ is nonincreasing with  $0\leq \widetilde{\theta}_0 \leq \widetilde{f}_0,$ and
	\begin{itemize}		\item[(e)] $\gamma^{\prime}>0 $ or $\widetilde\gamma >0,$ then  $|\nabla v| \in M^{\Phi , \zeta}_{\rm{loc}}(\Omega)$ for all $0\leq  \zeta\leq N-1;$
		\item[(f)]  $\gamma^{\prime} = \widetilde\gamma=0,$ then $|\nabla v| \in M^{\Phi , \zeta}_{\rm{loc}}(\Omega)$ for all $0\leq \zeta <  N.$
		\end{itemize}
\end{theorem}
{The proof of Theorem \ref{morreyg} is based on ideas introduced in \cite[Proposition 4.3]{queirozt}. The main obstacles to overcome regarding the proof of  Theorem \ref{morreyg} arise from treating the nonhomogeneous setting imposed by the operator considered, as well as the nonlinearities involved in the problem. These  mathematical difficulties differ from those encountered in \cite{queirozt}, since the authors dealt with linear operators, and the nonlinearities considered there are a particular case of those in Theorem \ref{morreyg}.}

 {As a consequence of Theorem \ref{morreyg} we have the next corollary, whose proof consists in an  application of the uniform $L^{\infty}$ estimates given by Theorem \ref{unifbound}, a careful adaptation to a nonhomogeneous setting  of the proof of \cite[Proposição 1.15]{tavares}. }

\begin{corollary}\label{holderg}
	Suppose that \eqref{condphi}, \eqref{hipF}, and \eqref{generalhyp} are in force with   in \eqref{generalhyp} that $h \in L^{\infty}(\Gamma).$  The assertions below hold for a minimizer $v$ of the functional \eqref{func}. If $y = Y^{\prime} $ is bounded or nondecreasing, and
		\begin{itemize}	
	\item[(a)]  $\gamma^{\prime}>0 $ or $\widetilde\gamma >0,$  then	$ v \in C^{0 , \sigma}_{\rm loc}(\Omega)$ for all $0< \sigma\leq \frac{\delta_0}{1 + \delta_0};$ 
	\item[(b)]   $\gamma^{\prime} = \widetilde\gamma=0,$ then	 $ v \in C^{0 , \sigma}_{\rm loc}(\Omega)$ for all $0<\sigma < 1.$
		\end{itemize}		
	On the other hand, if $y = Y^{\prime} $ is nonincreasing, with $-1< \widetilde{\theta}_0 \leq \widetilde{f}_0 <0$ or $-1 < \theta_0 < 0 \leq \widetilde{f}_0$, and
	\begin{itemize}		\item[(c)] $\gamma^{\prime}>0 $ or $\widetilde\gamma >0,$ then	$ v \in C^{0 , \sigma}_{\rm loc}(\Omega)$ for all $0< \sigma\leq \frac{\delta_0}{1 + \delta_0};$ 
		\item[(d)] $\gamma^{\prime} = \widetilde\gamma=0,$  then $ v \in C^{0 , \sigma}_{\rm loc}(\Omega)$ for all $0<\sigma < 1.$
	\end{itemize}	
		Finally, if $y = Y^{\prime} $ is nonincreasing with  $0\leq \widetilde{\theta}_0 \leq \widetilde{f}_0,$ and
			\begin{itemize}	
					\item[(e)]   $\gamma^{\prime}>0 $ or $\widetilde\gamma >0,$ then  $ v \in C^{0 , \sigma}_{\rm loc}(\Omega)$ for all $0< \sigma\leq\frac{\delta_0}{1 + \delta_0};$ 
	\item[(f)]  $\gamma^{\prime} = \widetilde\gamma=0,$ then $ v \in C^{0 , \sigma}_{\rm loc}(\Omega)$ for all $0<\sigma < 1 .$
		\end{itemize}	
	
\end{corollary}
{An important feature of Corollary \ref{holderg} is the hypotheses assumed in the structural conditions of the nonhomogeneous operator, which  allow the inclusion of the $p$-Laplacian operator in the singular case, that is,  $1 < p < 2$.  In \cite{piment}, H\"{o}lder estimates are not available for this case.}

In what follows we provide a BMO regularity result, which, under additional conditions, improve some items of Corollary \ref{holderg}.

\begin{theorem}\label{bmo1}
	Suppose that \eqref{condphi}, \eqref{hipF},  \eqref{generalhyp}  are in force with $h \in L^{\infty}(\Gamma),$  that $\gamma^{\prime} = \widetilde \gamma = 0,$  and
 additionally, one of the  following conditions:
\begin{itemize}	\item[$(\widetilde{a})$]  The function $y=Y'$ given in \eqref{hipF} is bounded;
	\item[$(\widetilde{b})$] $y$ is nonincreasing with $0\leq \widetilde{\theta}_0 \leq \widetilde{f}_0 ;$
	\item[$(\widetilde{c})$] $y$ is nondecreasing.
	\end{itemize}
Assume also that \eqref{ed} holds true, or that
$\Phi^{\prime} = \phi$ is convex with $\delta_0 >1$, where $\delta_0$ is given in \eqref{condphi}. 
Then $\tfrac{\partial v}{\partial x_i} \in {\rm{BMO}}_{\rm loc}(\Omega ), i=1,\cdots,N$. In particular $v$ is locally Log-Lipschitz continuous, that is , for any domain
	$\Omega^{\prime} \subset \subset\Omega$, there exists a constant $C>0$ such that  
	\begin{align}
		|v(x)-v(y)|\leq C|x-y||\log|x-y||,\notag
	\end{align}
	for all $x,y\in \Omega^{\prime}$. Therefore,
	$v\in C^{0,\sigma}_{\rm{loc}}(\Omega)$ for any $0<\sigma<1$. 
\end{theorem}
{
	Theorem \ref{bmo1} is based on the ideas of \cite{LQT,piment,ZZZ,ZFZ,zhengt}. Comparing to \cite{piment}, in our framework,  BMO estimates are derived for the minimizers of our functional  under more general hypotheses. 
	Moreover, as a consequence, local  Log-Lipschitz regularity and $C^{0,\sigma}_{\rm{loc}}$ for all $0<\sigma<1 $ estimates were proved. In view of the generality of the nonlinear  terms considered in \eqref{func}, the proof of Theorem \ref{bmo1} requires several distinct estimates compared to \cite{LQT,piment, ZZZ, ZFZ, zhengt}. For instance, since such a result is an application of  Theorem \ref{unifbound}, the boundedness of $h$ in $ L^{\infty} $ is required. We also point out that Theorem \ref{bmo1} complements \cite[Theorem 2, Corollary 1]{piment}.
}

Considering an unbounded interface data, it will be proved  Campanato and H\"older estimates. In this case, we assume  $h \in W^{1,p^{\prime} + \varepsilon}(\Omega)$, with $p^{\prime} + \varepsilon <N$, because if $p^{\prime} + \varepsilon \geq N$, from the usual embeddings, one has that $h \in C({\Omega}),$ which implies that $h\in L^\infty(\Gamma)$, as considered in the previous cases.

{
In what follows, we state Theorem \ref{last1}, which complements \cite[Theorem 3]{piment}, and establishes Campanato estimates for minimizers of the functional \eqref{func}. }

	\begin{theorem}\label{last1}Suppose $\delta_0 \geq 1$ and that \eqref{ed} holds with $2<p<N.$  Consider that $h \in W^{1,p^{\prime} + \varepsilon}(\Omega)$ with   $p^{\prime}  \leq   \varepsilon <  N - p^{\prime},$ $y = Y^{\prime}, f = F^{\prime}, \overline{f} = {\overline{F}}^{\prime}$ are bounded, and that $\lambda^{\prime} = \beta^{\prime} = \gamma^{\prime} = \widetilde\lambda = \widetilde\beta = \widetilde\gamma =0.$ Then, $\frac{\partial v}{\partial x_i} \in L^{1, N - N \frac{\varepsilon}{p^{\prime} + \varepsilon}}_{C}(\Omega^{\prime})$ for every domain $\Omega^{\prime} \subset \subset \Omega,$ and for all  $i=1,\cdots,N.$ Additionally,  if $\phi(t)= t^{p-1}$ and $\max \left\{p^{\prime} , \frac{N}{p-1}  - p^{\prime}\right\} <   \varepsilon<  N - p^{\prime}$, then  $\frac{\partial v}{\partial x_i}$ belongs to the Campanato space $ L^{p, N - N \frac{\varepsilon}{p^{\prime} + \varepsilon}}_{C,\rm{loc}}(\Omega), i =1,\cdots,N $, and consequently, $u \in C_{\rm loc}^{0,\alpha}(\Omega)$, where $\alpha = 1-\frac{N}{p+ \varepsilon(p-1)}$.
\end{theorem}

{Since the uniform estimates provided by Theorem \ref{unifbound} are no longer available in this setting, in order to prove Theorem \ref{last1}, several additional difficulties arise in our arguments. To overcome these challenges, we adapted the ideas from the proof of \cite[Theorem 4]{piment} to our framework. However, because we deal with a more general functional, our approach requires new estimates that do not appear in the proof of the aforementioned result.
	}

Regarding the applicability of the results, note that, as pointed out in \cite{zhengt}, the functions below satisfy the structural condition \eqref{condphi}.
\begin{itemize}
	
	\item[(i)]  $\phi(t) = t^{p-1}$ satisfies \eqref{condphi} with $\delta_0 =p-1$ and $g_0 =p-1$;
	
	\item[(ii)] $\phi(t)=(1+t)\ln(1+t)-t$ satisfies \eqref{condphi} with $\delta_0 =1$ and $g_0 =2$;
	
	\item[(iii)] $\phi(t)=\ln(1+at)+bt$ satisfies \eqref{condphi} with $\delta_0 = b$ and $g_0 = a+b$, where $a>0$, $b>0$;
	
	\item[(iv)] $\phi(t)=t^a \log_c(bt+d)$ satisfies \eqref{condphi} with $\delta_0 =a$ and $g_0 = a+ \frac{1}{\ln d}$, where $a,b>0$, $c,d>1$;
	
	\item[(v)] $\phi(t)= \frac{t}{\ln d}$ satisfies \eqref{condphi} with $\delta_0 =a- \frac{1}{\ln d}$ and $g_0 =a$, where $b>0$, $c,d>1$, $a>1$;
	
	\item[(vi)]
	\[
	\phi(t)= \begin{cases} 
		\log_c (bt+d)^p, & 0 \leq t < t_0 \\
		bt^q + c, & t \geq t_0
	\end{cases}
	\]
	satisfies \eqref{condphi} with $\delta_0 = \min\{p,q\} > 0$ and $g_0 = \max\{p,q\}$, where $a,b,c,p,q,t_0 > 0$ such that ${at^p = bt^q + c}$ and $apt^{p-1} = bqt^{q-1}$.
\end{itemize}
From the previous functions, it is clear how to obtain examples of nonlinearities $F,\bar F, Y$ in the definitions of functional in \eqref{func}, replacing $\delta_0$ by $1+\theta_0, 1 + \bar \theta_0, 1+ \widetilde \theta_0$, and $g_0$ by $1+f_0, 1 + \bar f_0, 1+ \widetilde f_0$, respectively.

Now, we point out the main contributions of this manuscript, along with the  mathematical difficulties that arose during the proofs of the  results described before.
First, we note that without considering a restriction on the sign of the functions $h$ and $\varphi$, the minimization problem \eqref{func} may admit a minimizer that changes its sign, which may give rise to a free boundary. Thus, \eqref{func} constitutes a two-phase free boundary problem, which has not been previously addressed in the literature under the hypotheses considered in this work. Moreover, the results obtained in this manuscript constitute initial steps in understanding the analytic and geometric properties of the free boundary associated to  the minimizers of the energy functional in \eqref{func}.

 The general hypotheses considered in the structural conditions of the operator that drives problem \eqref{func} allows accommodating a wide range of nonhomogeneous operators, including the $p$-Laplacian case. Additionally, the hypothesis in \eqref{hipF} incorporates  several nonlinearities with subcritical growth. Furthermore, we also consider  nonhomogeneous boundary conditions. These general hypotheses enable the inclusion of various boundary and interface conditions, producing a framework that encompasses  a diverse class of functionals, in particular, the one considered in \cite{piment}. The assumptions in this manuscript engender several mathematical obstacles in our study that  not arose for example in the references \cite{LQT,ZZZ,piment,queirozt,tavares,ZFZ,zhengt}.

The remainder of the paper is structured as follows: In Section 2, we provide some preliminary facts that will be frequently used throughout this work, including basic results regarding Campanato and BMO spaces, iteration, and preliminary facts involving Orlicz and Orlicz-Sobolev spaces. Additionally, we present specific results concerning to $\Phi$-Harmonic functions and technical estimates related to the nonlinearities under consideration in problem \eqref{func}. In Section 3, we present the proof of Theorem \ref{existence-boundness}, and adding the assumption $h \in L^{\infty}(\Omega)$, the proof of Theorem \ref{unifbound}. In Section 4, we prove the Morrey-type results, given by Theorems  \ref{morreyg} and Corollary \ref{holderg}. Section 5 is devoted to prove the BMO regularity in Theorem \ref{bmo1}. And the final section is dedicated to prove the case where $h$ belongs only to a suitable Sobolev space, specifically, the proof of Theorem \ref{last1}.

\section{Preliminaries}
\hfill

This section is dedicated to presenting preliminary facts that will be frequently used throughout this work.

First, we recall the definition of Campanato and BMO spaces. Let $\Omega \subset \mathbb{R}^N$ be a bounded domain.

\begin{definition}
Consider  $1 \leq p < \infty$, and $\lambda \geq 0$.  
The Campanato space $L^{p,\lambda}_{C}(\Omega)$ is defined as  
$$
L^{p,\lambda}_{C}(\Omega) := \left\{ f \in L^p(\Omega) ; [f]_{L^{p,\lambda}(\Omega)} < \infty \right\},
$$
where the Campanato semi-norm is given by
$$
[f]_{L^{p,\lambda}_{C}(\Omega)} := \sup_{\substack{x \in \Omega, r>0}} \left( \frac{1}{r^\lambda} \int_{\Omega \cap B_r(x)} |f(y) - f_{\Omega \cap B_r(x)}|^p \, dy \right)^{\frac{1}{p}},
\mbox{ with } f_{\Omega \cap B_r(x)} := \displaystyle\fint_{\Omega \cap B_r(x)}   f \text{d}x.$$
 It will be denoted by $L^{p,\lambda}_{C, \rm{loc}}(\Omega)$ the set of the functions which belongs to $L^{p,\lambda}_{C}(\Omega^{\prime})$ for all domain $\Omega^{\prime} \subset \subset \Omega. $
\end{definition}

\begin{definition}
	
	We denote by  $	\rm{BMO}(\Omega) $ the space
\begin{equation*}
	{\rm{BMO}}(\Omega) := \left\{ f \in L^1(\Omega) ; \sup_{\substack{x \in \Omega, r>0}} \fint_{\Omega \cap B_r(x)} |f -f_{\Omega \cap B_r(x)} | \text{d}x < \infty \right\}.
\end{equation*}
 It will be denoted by $\rm{BMO}_{\rm{loc}}(\Omega)$ the set of the functions which belongs to $\rm{BMO}(\Omega^{\prime})$ for all domain $\Omega^{\prime} \subset \subset \Omega. $
	
\end{definition}

The following result will play an important role in obtaining $L^{\infty}$ estimates for the minimizers of \eqref{func}.

\begin{lemma}\cite[Lemma 4.7]{LU}\label{lady-aux}
	Let $(P_n)$ be a sequence of non-negative numbers such that
	$$P_{n+1} \leq CD^{n}P^{1+\xi}_{n}, n\in \mathbb{N}\cup\{0\}, $$
	with $C,\xi>0$ and {$ D>1$} being constants that does not depend on $n.$ If
	$ P_0 \leq C^{- \frac{1}{\xi}} D^{-\frac{1}{\xi^2}},$
	then $P_n \rightarrow 0$ as $n \rightarrow \infty.$
\end{lemma}

Below we provide the particular case of \cite[Lemma 7.3]{giu} that will be used to obtain Morrey type estimates for the minimizers of \eqref{func}.

\begin{lemma}\label{glemma}
	Let $\psi : (0, R_0 ] \to \mathbb{R}$ be a non-negative and nondecreasing function. Suppose that there is  $0 < \tau < 1$ and constants $0<\beta < \delta$ and $B > 0$ such that
	$$
	\psi(\tau R) \leq \tau^\delta \psi(R) + B R^\beta,
	$$
	for each $0 < R < R_0$. Then there is a constant $C>0$, which depends on $\tau$, $\beta$, and $\delta,$ such that
	\[
	\psi(r) \leq C \left(\frac{r}{R}\right)^\beta \psi(R) + C r^\beta,
	\]
	for all $0 < r < R \leq R_0$.
\end{lemma}

Next result can be found in \cite[Lemma A.3]{zhengt}, it will be applied to obtain the local H\"{o}lder regularity to the minimizers of \eqref{func}.

\begin{lemma}\label{itg}
	\ \  {Let $\psi$ be a non-negative function on an interval $(0,R_0]$ with $R_0 \leq 1$. Suppose that
		\begin{align}
			\psi(r)\leq
			A\bigg(\frac{r}{R}\bigg)^{\alpha}\psi(R)+BR^{\beta},\notag
		\end{align}
		for all $0<r\leq R\leq R_{0}$, where  $A>1$ is a constant and $B, \alpha$ and $\beta$ are positive constants with $0<\beta<\alpha$. Fix $\delta \in (\beta,\alpha)$ and consider $\theta \in (0,1)$ with $A\theta^{\alpha} = \theta^{\delta}.$ Suppose that there exists a constant $m>0$ such that $\psi(r)\leq  m  \psi (\theta^{k}R)$ for all non-negative integer  $k$  and $r \in[\theta^{k+1} R,  \theta^{k} R].$ Then for any  $0<\sigma\leq \beta,$ there is a positive constant $C$, which depends on  $\phi, \theta$ and $m$ such that
		\begin{align}
			\psi(r) \leq  C\bigg(\frac{r}{R}\bigg)^{\sigma}\left(\psi(R)+BR^{\sigma}\right), \notag
		\end{align}
		for all $0< r \leq R \leq R_0$, and there is a positive constant $D>0,$which depends on $\phi, C,  R_0, \theta$ and $B$ such that
		\begin{align}
			\psi(r) \leq Dr^{\sigma} , \quad \mbox{	for all } 0<r\leq R_0. \notag
		\end{align}
	}
\end{lemma}

Now, several properties involving the nonlinearities of \eqref{func} are presented.
\begin{lemma}[\cite{zhengt}]\label{f-properties}{The functions $F$ and $f$ satisfy the following
		properties:}
	\begin{enumerate}
		\item[$(F_{1})$]
		${\min\{s^{1+\theta_0}, s^{1+f_0} \} F(t)\leq  F(st) \leq
			\max\{s^{1+\theta_0}, s^{1+f_0}\}F(t)},$ for all $s,t\geq 0;$
		\item[$(F_{2})$]${F(s+t)\leq 2^{1+f_0}(F(s)+F(t))},$ for all $s,t\geq 0;$
		\item[$(F_{3})$] If $f(t)$ is nonincreasing in $t>0$, then $F(s)-F(t)\leq F(s-t),$ for all $s,t \geq 0;$
		\item[$(F_{4})$] If $f(t)$ is nondecreasing in $t>0$, then $F(s)-F(t)\leq f(M)(s-t),$ for all $0\leq t\leq s\leq M,$ for some $M>0$;
		\item[$(f_{1})$]${\frac{1+\theta_0} {1+f_0}\min\{s^{\theta_0}, s^{f_0}\}f(t)\leq f(st)\leq \frac{1+f_0}{1+\theta_0}\max\{s^{\theta_0}, s^{f_0}\}f(t)},$ for all $s,t\geq 0;$
		\item[$(f_2)$]
		{If $f_0\leq \delta_0$, then $ f(t)\leq \frac{1+f_0}{1+\theta_0}\frac{f(1)}{g(1)}g(t),$ for all $  t\geq 1$.
		}
	\end{enumerate}
\end{lemma}

\begin{remark}\label{similar}
	Note that since \eqref{hipF} holds, then	 similar versions of Lemma \ref{f-properties} are valid  for $\overline{F},$ and $Y$ given in \eqref{hipF}.
\end{remark}

The proof of the  following result is included in    \cite[Theorem 1.1]{LQT}.
\begin{lemma}[\textnormal{\cite[pages 16,17]{LQT}}] \label{ineqpart} Consider $f_1, f_2 : \Omega \rightarrow \mathbb{R}$ measurable functions. It holds that $|f_{1}^{\pm} -f_{2}^{\pm}| \leq   |f_1 - f_2|.$ 
\end{lemma}

\subsection{Orlicz and Orlicz-Sobolev spaces}\label{suborlicz}
\hfill

Consider $\Omega \subset \mathbb{R}^{N}$ a bounded domain with smooth boundary, and the functions $\phi, \Phi$ given in Section \ref{intro}. 
The Orlicz space $L^{\Phi}(\Omega)$ is a Banach space defined as  the linear hull of the Orlicz class $$L_{\Phi}:= \{v: \Omega \rightarrow \mathbb{R} \ \text{measurable}; \int_{\Omega} \Phi(|v|)  \text{d}x < \infty\},$$
which is convex, due to the fact that $\Phi$  is convex with the norm
$$ \| v\|_{L^{\Phi}}:= \inf \left\{ \lambda > 0 ; \int_{\Omega} \Phi \left( \frac{|v|}{\lambda}\right) \leq 1 \right\}.$$
	
The Orlicz-Sobolev space $W^{1,\Phi}(\Omega)$ is  defined as the set of weakly differentiable functions $u \in W^{1,1}(\Omega)$ satisfying
$$ \int_{\Omega} \Phi (| v|) \ \text{d}x + \int_{\Omega} \Phi (|\nabla v|) \ \text{d}x < \infty,$$	
which becomes a Banach space with the norm $$\| v\|_{W^{1,\Phi}(\Omega)}:= \| v\|_{L^{\Phi}(\Omega)} + \| \nabla v\|_{L^{\Phi}(\Omega)}.$$  The space $W^{1,\Phi}_{0}(\Omega)$ is defined as the closure, with respect to the previous norm, of $C^{\infty}_{0}(\Omega).$
 Since  $\phi$ is strictly increasing there exists  its inverse function $\phi^{-1}$, hence the following Young-type inequality holds. Furthermore, in what follows, some useful results about Orlicz-Sobolev spaces are stated as matter of completeness.

\begin{lemma}[\cite{MW}]\label{tildes} The function $\phi$ admits it inverse $\phi^{-1},$ which  satisfies 
\begin{itemize}
	\item[(i)] 
	\begin{equation*}
	 ab \leq \varepsilon \Phi(a)+C(\varepsilon)\widetilde{\Phi}(b),
	\end{equation*}
for all  $a,b,\varepsilon >0,$ where  $\widetilde{\Phi}'(t)=\phi^{-1}(t),$ and $C(\varepsilon)$ is a constant which depends only on $\varepsilon.$
\item[(ii)] $\widetilde{\Phi}(\phi(t))\leq g_0 \Phi(t)$ for all $t>0.$
\end{itemize}	
\end{lemma}
\begin{lemma}[\cite{MW}]\label{modular} There exists a constant $C=C({\delta_0},g_0)$ such that
	\begin{align*}\| v\|_{L^{\Phi}(\Omega)}\leq C \max \left\{\left(\int_{\Omega} \Phi(|v|)\text{d}x \right)^{\frac{1}{1+{\delta_0}}}, \left(\int_{\Omega} \Phi(|v|)\text{d}x\right)^{\frac{1}{1+g_0}}  \right\}. \end{align*}
\end{lemma}

\begin{theorem}[\cite{MW}]  The dual of the space Orlicz space $L^{\Phi}(\Omega)$ is $L^{\widetilde{\Phi}}(\Omega),$ and  $L^{\Phi}(\Omega)$ and $W^{1,\Phi}(\Omega)$ are reflexive.
\end{theorem}

\begin{theorem}[\cite{MW}]\label{embedding}
It holds the continuous embedding 	$L^{\Phi}(\Omega) \hookrightarrow L^{1+{\delta_0}}(\Omega)$.	
\end{theorem}

\begin{lemma}[\cite{G}]\label{poincare} It  holds that $\int_{\Omega} \Phi(|v|) \text{d}x \leq  \int_{\Omega} \Phi(C|\nabla v|) \text{d}x, \forall v \in W^{1,\Phi}_{0}(\Omega)$, where the constant $C$ is twice the diameter of $\Omega.$
\end{lemma}

\begin{lemma}[\cite{MW}]\label{G-properties}The functions $\phi$ and $\Phi$ satisfy the following
	properties:
	\begin{enumerate}
		\item[$(\phi_{1})$] $\min\{s^{{\delta_0}}, s^{g_{0}} \}\phi(t) \leq  \phi(st) \leq
		\max\{s^{{\delta_0}}, s^{g_{0}}\}\phi(t),\ \text{for all} \ s,t\geq0$;
		\item[$(\Phi_{1})$] $\Phi \text{\ is\  convex on } [0,+\infty)\text{\ and }\ C^{2}-\text{continuous on}\  (0,+\infty)$;
		\item[$(\Phi_{2})$]$\frac{t\phi(t)}{1+g_{0}}\leq \Phi(t)\leq t\phi(t), \ \text{for all} \ t\geq
		0$;
		\item[$(\Phi_{3})$]$\min\{s^{{\delta_0}+1}, s^{g_{0}+1} \}\tfrac{\Phi(t)}{1+g_{0}}
		\leq \Phi(st) \leq
		(1+g_{0})\max\{s^{{\delta_0}+1}, s^{g_{0}+1}\}\Phi(t),\ \text{for all} \ s,t\geq 0.$
		\item[$(\Phi_{4})$] $\Phi(a+b)\leq 2^{g_{0}}(1+g_{0})(\Phi(a)+\Phi(b)), \ \text{for all} \
		a,b\geq
		0$;
	
	\end{enumerate}
\end{lemma}

\subsection{Results for $\Phi-$harmonic functions}
\hfill

This subsection describes some properties of  $\Phi-$harmonic functions that will be often used in this manuscript. 
Given $v \in W^{1,\Phi}(B_R),$  where $B_R := B_R (x_0)$ denotes the  open ball with radius $R$ and center $x_0$,  let $v^* \in W_v^{1,\Phi}(B_R)$ be the $\Phi-$harmonic replacement of  $v$ in $B_R(x_0),$ that is,   the unique minimizer weak solution of the problem

\begin{equation}\label{hreplace}
\left\{\begin{array}{rcl}
	\text{div} \ \dfrac{\phi(|\nabla v^*|)}{|\nabla v^*|}\nabla v^*&=&0 \quad \text{in }  B_R,\\
	v^*&=&v\quad \text{on } \partial B_R.
\end{array}
\right. 
\end{equation}

\begin{theorem}[\textnormal{\cite[Theorem 2.3]{MW}}] \label{simon}Let  $v,v^*\in W^{1,\Phi}(B_R)$ be as in \eqref{hreplace}.
Then

\begin{equation}\label{mart}
\int_{B_R} \Phi (|\nabla v|) - \Phi (|\nabla v^*|) dx \geq  C_3 \left( \int_{A_2} \Phi(|\nabla v - \nabla v^*|) \text{d}x + \int_{A_1} \dfrac{\phi(|\nabla v|)}{|\nabla v|} |\nabla v - \nabla v^* |^{2} \text{d}x\right),
\end{equation}
where 
$A_1 := \{x\in B_R ;  |\nabla v -\nabla v^*| \leq 2 |\nabla v|\},$ and $A_2 := \{x\in B_R ;  |\nabla v -\nabla v^*| > 2 |\nabla v|\},$ with the constants $C_i>0 , i=1,2,3$ not  depending on $R.$
\end{theorem}

The first assertion of the next result can be found in \cite[Lemma 2.5]{dasilva}, and the second one is contained in  in \cite[Page 6]{piment}, and  \cite[Lemma 3]{piment}. %\textcolor{red}{É a mesma referência? Devemos juntar?}

\begin{lemma}\label{simon2}
Consider   $v,v^*\in W^{1,\Phi}(B_R)$  as in \eqref{hreplace}.
The following assertions hold.
\begin{itemize}\item[(i)] If $\Phi^{\prime} = \phi$ is convex,  there is a constant $C>0$, which depends only on $\delta_0>1, $ and $g_0$, such that
\begin{equation*}
\int_{B_{R}} \Phi (|\nabla v - \nabla v^*|) \textnormal{d}x \leq C \int_{B_{R}} \Phi(|\nabla v|) - \Phi(|\nabla v^*|) \textnormal{d}x.
\end{equation*}	
\item[(ii)] If \eqref{ed} holds, then $v,v^* \in W^{1,p}(B_R(x_0))$, and  there is a constant $C>0$, which does not depend  on $R$, such that
\begin{equation*}
	\int_{B_R} |\nabla v - \nabla v^*|^{p} \textnormal{d}x\leq C \int_{B_R} \Phi(|\nabla v|) - \Phi(|\nabla v^*|) \textnormal{d}x
\end{equation*}	

\end{itemize}

\end{lemma}
{
\begin{lemma}[\textnormal{\cite[Page 345]{L}}]\label{iotal0}
	Consider $v,v^* \in W^{1,\Phi}(B_R)$ as in \eqref{hreplace},  and  suppose that   ${v^* \in L^{\infty}(B_R(x_0))}$ with   $\| v^*\|_{L^{\infty}(B_R)} \leq K$.
 Thus, it follows that
	\begin{equation}\label{hreplacement}
		\int_{B_r} \Phi (|\nabla v^*|)\text{d}x \leq C \left( \frac{r}{R}\right)^N \int_{B_R }\Phi (|\nabla v^*|)\text{d}x, \ \mbox{for} \ 0<r \leq R,
	\end{equation}
	where $C=C(N,\delta_0,g_0,K)$ is a positive constant.
	\end{lemma}}

\begin{lemma}[\textnormal{\cite[(5.9), page 346]{L}}]\label{iotal}
	Consider $v,v^* \in W^{1,\Phi}(B_R)$ as in \eqref{hreplace}, and  suppose that   ${v^* \in L^{\infty}(B_R(x_0))}$ with   $\| v^*\|_{L^{\infty}(B_R)} \leq K$. Then, for every $\mu \in (0,N)$, there exists a constant ${C=C(\mu,N,\delta_0,g_0,K)>0}$ such that
	\begin{align*}
	\int_{B_{r}}\Phi(|\nabla v^*|) \textnormal{d}x\leq Cr^{\mu}, \quad \mbox{ for all } \ 0<r \leq R.
	\end{align*}
\end{lemma}

The result below will play an important role, and its proof will be provided for convenience.
\begin{lemma}\textnormal{\cite[Lemma 3.1]{ZFZ}}\label{liota}
	Consider $v,v^* \in W^{1,\Phi}(B_R)$ as in \eqref{hreplace}, and  suppose that   $v^* \in L^{\infty}(B_R)$ with   $\| v^*\|_{L^{\infty}(B_R)} \leq K$. Then for any $\mu \in (0,N)$, there exists a constant  $C=C(\mu,N,{\delta_0},g_0,K)>0$ such that
	\begin{align*}
	\int_{B_R}\Phi(|\nabla v-\nabla v^*|)\textnormal{d}x\leq&  C\int_{B_R}(\Phi(|\nabla v|)-\Phi(|\nabla v^*|))\textnormal{d}x\notag +CR^{\frac{\mu}{2}}\bigg(\int_{B_R}(\Phi(|\nabla v|)-\Phi(|\nabla v^*|))\textnormal{d}x\bigg)^{\frac{1}{2}}.
	\end{align*}

\end{lemma}
\begin{proof}
	
Consider $A_1 ,$ and $A_2$ as in Theorem \ref{simon}. From \eqref{mart}, we obtain for some positive constant $C>0$, that

\begin{equation}\label{A1}
C\int_{A_1}\dfrac{\phi(|\nabla v|)}{|\nabla v|} |\nabla v - \nabla v^* |^{2} \text{d}x \leq \int_{B_R} \Phi(|\nabla v|) - \Phi(|\nabla v^*|)\text{d}x,
\end{equation}
and
\begin{equation*}
C\int_{A_2} \Phi(|\nabla v - \nabla v^*|)\text{d}x \leq \int_{B_R} \Phi(|\nabla v|) - \Phi(|\nabla v^*|)\text{d}x.
\end{equation*}

In addition, we have
\begin{equation*}
	\int_{B_R}  \Phi (|\nabla v - \nabla v^*|) \text{d}x = \int_{A_1 } \Phi (|\nabla v - \nabla v^*|) \text{d}x  + \int_{A_2 } \Phi (|\nabla v - \nabla v^*|) \text{d}x .
\end{equation*}

%Note that $v^*$ is the unique minimizer of the energy functional $J(w):= \int_{B_R}\Phi(|\nabla w|)\text{d}x $ defined over the set $\mathcal{K}:=\{w \in W^{1,\Phi} (B_R (x_0)) ; w - v \in W^{1,\Phi}_{0} (B_R (x_0)) \}.$ Thus, \begin{equation}\label{diri}\int_{B_R} \Phi(\nabla w) - \Phi (\nabla v^*)\text{d}x \geq 0,\end{equation}for all $w\in \mathcal{K}.$ 

From $(\Phi_2)$ of Lemma  \ref{G-properties}, it follows that $\left(\tfrac{\Phi(t)}{t}\right)^{\prime}\geq 0,$ for all $t>0$. Thus, it yields $\tfrac{\Phi(t)}{t}$ is nondecreasing for  $t>0$, which combined with $(\Phi_2),(\Phi_3),$ {the Cauchy-Schwarz inequality, and \eqref{A1} provides}
\begin{small}
\begin{equation*}
\begin{aligned}
\int_{A_1} \Phi (|\nabla v - \nabla v^*|)\text{d}x  &= \int_{A_1} \frac{\Phi(|\nabla v - \nabla v^*|)}{|\nabla v - \nabla v^*|}|\nabla v - \nabla v^*|\text{d}x \leq \int_{A_1} \frac{\Phi(2|\nabla v|)}{2|\nabla v|}|\nabla v - \nabla v^*|\text{d}x\\
& \leq C  \int_{A_1} \dfrac{\Phi(|\nabla v|)}{|\nabla v|} |\nabla v - \nabla v^*|\text{d}x = \int_{A_1} \dfrac{\sqrt{\Phi(|\nabla v|)}}{|\nabla v|} |\nabla v - \nabla v^*| \sqrt{\Phi(|\nabla v|)} \text{d}x\\
& \leq C \left( \int_{A_1} \dfrac{\Phi(|\nabla v|)}{|\nabla v|^{2}} |\nabla v - \nabla v^*|^2\text{d}x\right)^{\frac{1}{2}} \left( \int_{A_1} \Phi (|\nabla v|)\text{d}x \right)^{\frac{1}{2}}\\
& \leq C \left(\int_{A_1} \dfrac{\phi(|\nabla v|)}{|\nabla v|} |\nabla v - \nabla v^*|^2\text{d}x\right)^{\frac{1}{2}}\left(\int_{A_1} \Phi(|\nabla v|)\text{d}x\right)^{\frac{1}{2}}\\
& \leq C \left( \int_{B_R} \Phi (|\nabla v|) - \Phi(|\nabla v^*|)\text{d}x\right)^{\frac{1}{2}}\left(\int_{B_R} \Phi(|\nabla v|)\right)^{\frac{1}{2}}\\
& = C \left( \int_{B_R} \Phi (|\nabla v|) - \Phi(|\nabla v^*|)\text{d}x \right)^{\frac{1}{2}}\left(\int_{B_R} \Phi (|\nabla v|) - \Phi (|\nabla v^*|) + \Phi(|\nabla v^*|) \text{d}x\right)^{\frac{1}{2}}\\
& \leq C \left( \int_{B_R} \Phi (|\nabla v|) - \Phi(|\nabla v^*|)\text{d}x \right)^{\frac{1}{2}} \left[\left( \int_{B_R} \Phi (|\nabla v|) - \Phi(|\nabla v^*|)\text{d}x \right)^{\frac{1}{2}} + \left(\int_{B_R} \Phi(|\nabla v^*| \text{d}x\right)^{\frac{1}{2}}\right]\\
&\leq C \left( \int_{B_R} \Phi (|\nabla v|) - \Phi(|\nabla v^*|) \text{d}x\right) + C \left( \int_{B_R} \Phi (|\nabla v|) - \Phi(|\nabla v^*|) \text{d}x\right)^{\frac{1}{2}} \left(\int_{B_R} \Phi(|\nabla v^*| \text{d}x\right)^{\frac{1}{2}}.
\end{aligned}
\end{equation*} 
\end{small}
The last inequality, and  Lemma \ref{iotal} provide the result.
	
\end{proof}

The  next result, whose first part can be found in \cite[Lemma 5.1]{L} and in \cite[Lemma A.16]{zhengt}, and the second one in \cite[Lemma 4.1]{baroni} provides an important property of $\Phi$-harmonic functions, and plays an essential role in the development of $C^{1,\sigma}_{loc}$ estimates for such functions.

\begin{lemma} \label{liebreg} Let $v^* \in W^{1,\Phi}(B_R)$  be a weak solution of \eqref{hreplace}.
The following assertions are true.
\begin{itemize}
\item[(a)] It holds that 	\begin{align*}	
		\int_{B_r}\Phi(|\nabla v^* - (\nabla v^*)_{x_0,r}|) \textnormal{d}x \leq C \left( \frac{r}{R}\right)^{{N}+\sigma} \int_{B_R}\Phi(|\nabla v^* - (\nabla v^*)_{x_0 ,R}|) \textnormal{d}x,
	\end{align*}
	for all $0<r\leq R,$ where $C>0$ is constant that depends only on $N,$ and $g_0.$ Moreover, $\sigma \in (0,1)$ is a constant which does not depend on $r$ and $v^*,$ and the symbol $(\nabla v^*)_{x_0 , r} $ is defined by
	$$ (\nabla v^*)_{x_0 , r} : = \left(\fint_{B_r (x_0)} \frac{\partial v^* }{\partial x_1} \textnormal{d}x , \cdots, \fint_{B_r (x_0)} \frac{\partial v^* }{\partial x_N} \textnormal{d}x \right).$$
\item[(b)]	 If  \eqref{ed} holds, then
$$		\int_{B_r(x_0)}|\nabla v^* - (\nabla v^*)_{x_0,r}| \textnormal{d}x \leq C \left( \frac{r}{R}\right)^{{N}+\sigma} \int_{B_R}|\nabla v^* - (\nabla v^*)_{x_0 ,R}| \textnormal{d}x,$$
for all $0 < r \leq R$, and some constant $C>0$ with $\sigma \in (0,1]$ not depending on $v^*$ and $r.$
\end{itemize}
\end{lemma}

The result below, whose  first part can be found for example in \cite[Lemma 3.1]{ZFZ} or \cite[Lemma 4.1]{ZZZ}, and the second one in \cite[Proposition 2]{piment}, for  $\Phi$-harmonic functions will play an important role in obtaining the H\"{o}lder regularity to the gradient of the minimizers of \eqref{func}.

	\begin{lemma} \label{comp} Consider $v,v^* \in W^{1,\Phi}(B_R)$ as in \eqref{hreplace}.
The following assertions are true.	
\begin{itemize}
\item[(a)]It holds that	
	\begin{equation*}
\int_{B_r} \Phi(|\nabla v-(\nabla v)_{x_0,r}|) \textnormal{d}x\leq
	C\bigg(\frac{r}{R}\bigg)^{N+\sigma}\int_{B_R} \Phi(|\nabla v-(\nabla
	v)_{x_0,R}|) \textnormal{d}x+C\int_{B_R} \Phi(|\nabla v-\nabla v^*|) \textnormal{d}x
\end{equation*}
	for all $0<r\leq R,$ where $C>0$ is constant that depends only on $N,$ and $g_0,$ with $\sigma \in (0,1)$ being the  constant given in Lemma \ref{liebreg}.
\item[(b)] If additionally  \eqref{ed} holds, then
	\begin{equation*}
	\int_{B_r} |\nabla v-(\nabla v)_{x_0,r}| \textnormal{d}x\leq
		C\bigg(\frac{r}{R}\bigg)^{N+\sigma}\int_{B_R} |\nabla v-(\nabla
		v)_{x_0,R}| \textnormal{d}x+C\int_{B_R} |\nabla v-\nabla v^*| \textnormal{d}x,
\end{equation*}
\end{itemize}
for all $0<r\leq R,$ where $  \sigma \in (0,1)$ was given in Lemma \ref{liebreg}, and  some constant $C>0$, which does not depend  on $v^*,$ and $R.$
\end{lemma}

\section{Existence of Minimizers}
\hfill

The goal of this section is to prove that the energy functional associated with \eqref{func} admits at least one minimizer.
Before establishing the existence of minimizers, we will first show in the result below that the functional \eqref{func} is well-defined.

\begin{proposition} The functional \eqref{func} is well defined.
\end{proposition}
\begin{proof} We claim  that $\displaystyle\int_{\Gamma} Y(v^{+}) \text{d}H^{N-1}$ is finite.  First, we observe that,  from Remark \ref{similar} it follows that 
	\[0\leq Y(t) \leq C (t^{1+\widetilde{\theta}_0} + t^{1 + \widetilde{f}_0}), \quad \mbox{for all }t\geq  0,\]
for a positive constant $C$.	Moreover, 
	 if  $0<1+\widetilde{\theta}_0 \leq 1 +\widetilde{f}_0 < 1 $, we note that for any $0<\rho<1$, it holds  $t^{\rho} \leq 1 + t,$ for all $t > 0.$ Otherwise, there exists $\rho' >1$ such that 
 $1\leq 1+\widetilde{\theta}_0 \leq 1 +\widetilde{f}_0 \leq \rho' < (1+\delta_0)^{*}$ then,  
\begin{equation}\label{boa1}
0\leq Y(t) \leq C (t^{1+\widetilde{\theta}_0} + t^{1 + \widetilde{f}_0}) \leq C(1 + t + t^{\rho^{\prime}}),
\end{equation}	
for all $t\geq  0$.
 Thus, from  the continuous embeddings  
  $W^{1,\Phi}(\Omega) \hookrightarrow W^{1,1+\delta_0}(\Omega) \hookrightarrow L^{\widetilde{\rho}}(\Gamma), 1 \leq \widetilde{\rho} \leq (1+\delta_0)_{*}$
  we have that $\int_{\Gamma} Y(v^{+})\text{d}H^{N-1}$ is finite.
 The previous reasoning  implies that 
 $$\displaystyle\int_{\Gamma} Y(v^{-}) \text{d}H^{N-1}, \displaystyle\int_{\Omega_1} F(v^{+}) \text{d}x \displaystyle\int_{\Omega_1} F(v^{-}) \text{d}x, \displaystyle\int_{\Omega_2} \overline{F}(v^{+}) \text{d}x, \mbox{ and } \displaystyle\int_{\Omega_2} \overline{F}(v^{-}) \text{d}x   \mbox{ are finite.} $$
 Furthermore, in view of the continuous embeddings $W^{1,1+\delta_0}(\Omega) \hookrightarrow W^{1,1}(\Omega),$ and $W^{1,1}(\Omega_i)\hookrightarrow L^{\frac{N}{N-1}}(\Omega_i),i=1,2$ and $g,\overline{g}\in L^{{N}}(\Omega_i),i=1,2,$ we obtain that $\int_{\Omega_1} g v \text{d}x, $ and  $\int_{\Omega_2} \overline{g} v \text{d}x $  are finite.  The finiteness of the remaining terms in the functional \eqref{func} is clear. Hence, the result is established.
\end{proof}

\begin{proof}[Proof of Theorem \ref{existence-boundness}]
We will begin the proof by showing the existence of a minimizer. We claim that 
\begin{equation}\label{J_0}
 J_0:=\displaystyle\inf_{v \in W^{1,\Phi}_{\varphi}(\Omega)}  J(v) > -\infty.
 \end{equation}
In fact, let $v \in W^{1,\Phi}(\Omega).$ Since   $v\in W^{1,1}(\Omega_1)\hookrightarrow L^{\frac{N}{N-1}}(\Omega_1), $ we obtain 	
\begin{equation*}%\label{eps1}
\begin{aligned}
\int_{\Omega_1} |g||v| \text{d}x  & \leq \|g\|_{L^N(\Omega_1)}\|v\|_{L^{1^*}(\Omega_1)}  \leq  \|g\|_{L^N(\Omega_1)}(\|v-\varphi\|_{L^{1^*}(\Omega_1)}+\|\varphi\|_{L^{1^*}(\Omega_1)})\\
& \leq  \|g\|_{L^N(\Omega_1)}(\|v-\varphi\|_{L^{1^*}(\Omega)}+\|\varphi\|_{L^{1^*}(\Omega)})  \leq \|g\|_{L^N(\Omega_1)}(C\|\nabla v-\nabla\varphi\|_{L^{1}(\Omega)}+\|\varphi\|_{L^{1^*}(\Omega)}) \\
& \leq C_1 \|g\|_{L^N(\Omega_1)}(\|\nabla v\|_{L^{1}(\Omega)} + \|\nabla\varphi\|_{L^{1}(\Omega)}+\|\varphi\|_{L^{1^*}(\Omega)}) \\
& \leq { C_2\|g\|_{L^N(\Omega_1)}\left(\varepsilon\int_{\Omega} \Phi(|\nabla v|) \text{d}x + C_{\varepsilon}\int_{\Omega} \widetilde{\Phi}(1) \text{d}x +\|\nabla\varphi\|_{L^{1}(\Omega)}+\|\varphi\|_{L^{1^*}(\Omega)}\right)},
\end{aligned}
\end{equation*}
for some constants $C_1, C_2, C_{\varepsilon},K>0,$ where in the last inequality we applied Lemma \ref{tildes}.
Applying the same reasoning for the term $\int_{\Omega_2} |\overline{g}||v| \text{d}x$ it follows that 
\begin{equation}\label{y1}
\int_{\Omega_1} |g||v| \text{d}x + \int_{\Omega_2} |\overline{g}||v| \text{d}x \leq { C(\|g\|_{L^N(\Omega_1)} + \|\overline{g}\|_{L^N(\Omega_2)})\left(\int_{\Omega}\varepsilon \Phi(|\nabla v|) + C_{\varepsilon} \widetilde{\Phi}(1) \text{d}x +\|\nabla\varphi\|_{L^{1}(\Omega)}+\|\varphi\|_{L^{1^*}(\Omega)}\right)},
\end{equation}
for some constants $C,C_{\varepsilon}>0.$
From \eqref{generalhyp}, and Lemma \ref{embedding} we have the continuous embeddings  
$$W^{1,\Phi}(\Omega_1) \hookrightarrow W^{1,1+\delta_0}(\Omega_1) \hookrightarrow L^{(1+\delta_0)_{*}}(\Gamma) \hookrightarrow L^{m}(\Gamma).$$
Thus, using Minkowski Inequality and  applying the Young Inequality  for $\tau = \frac{1}{1+\delta_0}$,  we obtain that
\begin{equation}\label{eps2}
\begin{aligned}
\int_{\Gamma} |h||u| \text{d}x & \leq \| h\|_{L^{\frac{m}{m-1}}(\Gamma)} \| v\|_{L^{m}(\Gamma)}  \leq C  \|v  \|_{L^{(1+\delta_0)_{*}}(\Gamma)} \leq  C\| v \|_{W^{1,1+\delta_0}(\Omega_1)} \\
& \leq C(\| v - \varphi\|_{W^{1,1+\delta_0}(\Omega_1)}  + \|\varphi \|_{W^{1,1+\delta_0}(\Omega_1)})\leq  C (\| v - \varphi\|_{W^{1,1+\delta_0}_{0}(\Omega)}  + \|\varphi \|_{W^{1,1+\delta_0}(\Omega)}) \\ &= C\left( \left( \int_{\Omega} |\nabla v - \nabla \varphi|^{1+\delta_0}\right)^{\frac{1}{1+\delta_0}} + \| \varphi\|_{W^{1,1+\delta_0}(\Omega)}\right) \leq C \left(\left( \int_{\Omega} |\nabla v|^{1+\delta_0}\right)^{\frac{1}{1+\delta_0}} + \| \varphi\|_{W^{1,1+\delta_0}(\Omega)}\right) \\ & \leq C\left( \tau \int_{\Omega} |\nabla v  |^{1+\delta_0} \text{d}x + K_{\tau} +\| \varphi\|_{W^{1,1+\delta_0}(\Omega)}\right) \leq C\left( \tau \int_{\Omega} \Phi(|\nabla v|) \text{d}x +C_{\tau}+\|\varphi \|_{W^{1,1+\delta_0}(\Omega)}\right),
\end{aligned}
\end{equation}
for positive constants $K_{\tau}, C_{\tau},C >0,$  where in the last inequality we used that
\begin{equation}\label{auxineq}
	|s|^{1+\delta_0} -1 \leq \frac{1+g_0}{\Phi(1)}\Phi(|s|), \ \text{for all} \ s\in\mathbb{R},
\end{equation}
which is a consequence of $(\Phi_3)$ in Lemma \ref{G-properties}.
Using the fact that $W(v), \overline{W}(v), Z(v) \geq 0,$ and the inequalities \eqref{y1}, \eqref{eps2}, it follows that \eqref{J_0} holds true. Then, we may consider a minimizing sequence $(v_j) $ and $j_0 \in \mathbb{N}$ satisfying $J(v_j) \leq J_0 + 1$ for all $j \geq j_0.$ In view of  \eqref{y1}, \eqref{eps2}, and using that $W(v), \overline{W}(v), Z(v)$ are non-negative functions, and considering a suitable choice of $\varepsilon,\tau  >0$ we obtain a $\rho>0$ small enough such that
\begin{equation}\label{refle1}
\begin{aligned}
		\int_{\Omega} \Phi(|\nabla v_j|) \text{d}x &= J(v_j) - \int_{\Omega_1}  W(v_j) + gv_j\text{d}x - \int_{\Omega_2}  \overline{W}({v_j}) + \overline{g}v_j\text{d}x + \int_{\Gamma} hv \ \text{d}H^{N-1} - \int_{\Gamma} Z(v) \ \text{d}H^{N-1}\\
	&	\leq \rho \int_{\Omega}  \Phi(|\nabla v_j|)+  C\left( \|\varphi \|_{W^{1,1+\delta_0}(\Omega)} + \|\nabla\varphi\|_{L^{1}(\Omega)}+  \|\varphi\|_{L^{1^*}(\Omega)} +1  \right) ,
\end{aligned}
\end{equation}
for some  constant $C >0,$ which does not depend on $j \in \mathbb{N}$.
Thus, by applying $(\Phi_3)$ and $(\Phi_4)$ of Lemma \ref{G-properties}, along with Lemma \ref{poincare} and the nondecreasing monotonicity of the function $\Phi$, it follows that
\begin{align}\label{refle2}
	\int_{\Omega} \Phi(|v_j|) \text{d}x &\leq   %\int_{\Omega} \Phi(|v_j-\varphi| + |\varphi|) \text{d}x\notag
	 C \left( \int_{\Omega} \Phi(|v_j-\varphi|) \text{d}x + \int_{\Omega} \Phi(|\varphi|) \text{d}x\right)\notag\\
	&\leq   C\left( \int_{\Omega}\Phi(|\nabla v_j - \nabla \varphi|) \text{d}x + \int_{\Omega} \Phi(|\varphi|) \text{d}x\right)\notag\\
&	\leq  C\left( \int_{\Omega}\Phi(|\nabla v_j|)\text{d}x + \int_{\Omega} \Phi(|\nabla \varphi|)  \text{d}x + \int_{\Omega} \Phi(|\varphi|) \text{d}x\right),
\end{align}
where $C>0$ is a constant depending  only on the diameter of $\Omega,$  ${\delta_0}$, and $g_0.$ Thus, from \eqref{refle1}, and \eqref{refle2} we obtain that  the sequence $(v_j)$ is bounded in $W^{1,\Phi}(\Omega).$
From $(\Phi_4)$ in Lemma \ref{G-properties}, Lemma \ref{poincare}, and  the nondecreasing monotonicity of  $\Phi$ we deduce that $v_j-\varphi$ is a bounded sequence in $W^{1,\Phi}_{0}(\Omega).$ Since the space $W^{1,\Phi}_{0}(\Omega)$ is reflexive, there is a function $v \in W^{1,\Phi}(\Omega)  $ such that $v-\varphi\in W^{1,\Phi}_0(\Omega)$, and  $v_j \rightharpoonup v \ \text{in} \ W^{1,\Phi}(\Omega)$, up to a subsequence.  
Fix $\max\{1, 1 + f_0, 1+\overline{f}_0\} <\eta < (1+\delta_0)^*,$ and let $\rho^{\prime}$ be as in \eqref{boa1}. From Lemma \ref{embedding}, we obtain that $v_j \rightharpoonup v \ \text{in} \ W^{1,1+\delta_0}(\Omega).$ Using the continuous embeddings $$W^{1,1+\delta_0}(\Omega_i)\hookrightarrow W^{1,1}(\Omega_i)\hookrightarrow L^{\frac{N}{N-1}}(\Omega_i), \, i=1,2,$$ and the  compact embeddings $W^{1,1+\delta_0}(\Omega) \hookrightarrow L^{r}(\Omega),$ for $r = \eta, \widetilde{q}, \rho^{\prime},$ up to a subsequence,  we conclude that
\begin{equation}\label{convers}
\begin{aligned}
	v_j &\rightarrow v \ \text{a.e. in} \ \Omega , \quad
	v_j \rightarrow v \ \text{a.e. in} \ \Gamma , \quad
	v_j  \rightharpoonup v \ \text{ in } \ L^{\frac{N}{N-1}}(\Omega_i),i=1,2,  \quad 	v_j \rightarrow v \ \text{ in } \ \ L^{\widetilde{q}}(\Gamma) ,\\
	|v_j|&\leq h_{\eta} \ \text{a.e. in } \ \Omega, \ \text{for some} \ h_{\eta} \in L^{\eta}(\Omega), \quad
		|v_j|\leq h_{{\rho}^{\prime}} \ \text{a.e. in } \ \Gamma, \ \text{for some} \ h_{{\rho}^{\prime}} \in L^{{\rho}^{\prime}}(\Gamma).
\end{aligned}
\end{equation}

{Since $\Phi$ is convex, and $|\nabla v_j| \rightharpoonup |\nabla v| \ \text{in} \ L^{\Phi}(\Omega),$  as a consequence of the Ioffe-De Giorgi lower semicontinuity theorem \cite[Theorem 1]{Ioffe1977}, it follows that 
\begin{align}\label{lim-inf}
	\int_{\Omega}\Phi(|\nabla v |) \text{d}x \le \liminf_{j \rightarrow \infty} \int_{\Omega}\Phi(|\nabla v_j |) \text{d}x.
\end{align}}
%In order to prove the claim,  the convexity of $\Phi$ provides that 
%\begin{align}\label{convexity}
%	\int_{\Omega} \Phi(|\nabla v_j|) \text{d}x \geq \int_{\Omega} \Phi(|\nabla v|)\text{d}x + \int_{\Omega} \phi(|\nabla v|)(|\nabla v_j| - |\nabla v|) \text{d}x.
%\end{align}
%Moreover, in view of Lemma \ref{tildes}, we know that$ %\widetilde{\Phi}\left(\phi(|\nabla v|) \left|\frac{\partial v }{\partial x_i} \frac{1}{|\nabla v|}\right|\right) \leq 	\widetilde{\Phi}(\phi(|\nabla v|)) \leq C \Phi(|\nabla v|),$which ensures  $\phi(|\nabla v |) %\dfrac{\nabla v}{|\nabla v|} \in L^{\widetilde{\Phi}}(\Omega).$ Since$|\nabla v_j| \rightharpoonup |\nabla v| \ \text{in} \ L^{\Phi}(\Omega),$ in view of the weak convergence,  passing \eqref{convexity} the limit, we obtain \eqref{lim-inf}.

 Now, we observe that, arguing as in \eqref{boa1},  we obtain the inequalities
 \begin{equation*}
 	\begin{aligned}
 	0\leq Y(t), F(t), \overline{F}(t)  \leq C(1 + t + t^{\eta}),
 	\end{aligned}
 \end{equation*}
  for all $t\geq 0,$ and a constant $C>0$. Hence, from \eqref{convers}, and the Lebesgue's Dominated Convergence Theorem it follows that
 \begin{equation}\label{Fatou1}
\begin{aligned}
	\lim_{j \rightarrow +\infty} \int_{\Omega_1} (F(v^{+}_{j})+F(v^{-}_{j}))\text{d}x &=  \int_{\Omega_1} (F(v^{+})+F(v^{-}))\ \text{d}x, \\
		\lim_{j \rightarrow +\infty} \int_{\Omega_2} (\overline{F}(v^{+}_{j})+\overline{F}(v^{-}_{j}))\text{d}x & =  \int_{\Omega_2}(\overline{F}(v^{+})+\overline{F}(v^{-}))\ \text{d}x,\\
			\lim_{j \rightarrow +\infty} \int_{\Gamma} (Y(v^{+}_{j}) + Y(v^{-}_{j}))\text{d}H^{N-1} & = \int_{\Gamma}(Y(v^{+}) + Y(v^{-}))\text{d}H^{N-1}.
\end{aligned}
\end{equation}
Note that 
\begin{equation*}
\begin{aligned}
\int_{\Omega_1}\chi_{\{v\leq 0\}} \text{d}x & =\int_{\Omega_1 \cap\{v\leq 0\}} \chi_{\{v_j >0\}} \text{d}x +  \int_{\Omega_1 \cap\{v\leq 0\}} \chi_{\{v_j \leq 0\}} \text{d}x \\
& \leq \int_{\Omega_1 \cap\{v\leq 0\}} \chi_{\{v_j >0\}} \text{d}x +  \int_{\Omega_1} \chi_{\{v_j \leq 0\}} \text{d}x.
\end{aligned}
\end{equation*}
Since $v_j \rightarrow v$ a.e. in $\Omega$,  from the Lebesgue's Dominated Convergence Theorem, we obtain 
\begin{equation}\label{inf0}
\int_{\Omega_1}  \chi_{\{v\leq0\}}\text{d}x \leq \lim_{j \rightarrow +\infty} \left [\int_{\Omega_1 \cap\{v\leq 0\}} \chi_{\{v_j >0\}} \text{d}x +  \int_{\Omega_1}\chi_{\{v_j \leq 0\}} \text{d}x\right],
\end{equation}
and also that 
\begin{equation*}
\begin{aligned}
\int_{\Omega_1}  \chi_{\{v>0\}} \text{d}x &
=\int_{\Omega_1 \cap\{v>0\}} \left( \displaystyle\lim_{j \rightarrow +\infty} \chi_{\{v_j >0\}}\right) \text{d}x = \lim_{j \rightarrow +\infty} \int_{\Omega_1 \cap\{v>0\}}  \chi_{\{v_j >0\}} \text{d}x,
\end{aligned}
\end{equation*}
which combined with \eqref{inf0} provides that
\begin{equation}\label{inf1}
\int_{\Omega_1} \lambda^{\prime}  \chi_{\{v>0\}} + \widetilde\lambda\chi_{\{v\leq0\}}\text{d}x \leq \lim_{j\rightarrow +\infty} \int_{\Omega_1} \lambda^{\prime}  \chi_{\{v_j>0\}} + \widetilde\lambda\chi_{\{v_j\leq0\}} \text{d}x.
\end{equation}
Analogously, we get that
\begin{equation}\label{inf2}
\begin{aligned}
	\int_{\Omega_2} \beta^{\prime}  \chi_{\{v>0\}} + \widetilde\beta\chi_{\{v\leq0\}} \text{d}x &\leq \lim_{j\rightarrow +\infty} \int_{\Omega_2} \beta^{\prime}  \chi_{\{v_j>0\}} + \widetilde\beta\chi_{\{v_j\leq0\}} \text{d}x,\\
	\int_{\Gamma} \gamma^{\prime} \chi_{\{v>0\}} + \widetilde\gamma\chi_{\{v\leq0\}} \text{d}x &\leq \lim_{j\rightarrow +\infty} \int_{\Gamma} \gamma^{\prime}  \chi_{\{v_j>0\}} + \widetilde\gamma\chi_{\{v_j\leq0\}} \text{d}x.
	\end{aligned}
\end{equation}
From  \eqref{convers} we obtain that 
\begin{equation*}
	\int_{\Omega_1}gv_j  \text{d}x \rightarrow  \int_{\Omega_1} g v \text{d}x, \int_{\Omega_2}\overline{g} v_j  \text{d}x \rightarrow  \int_{\Omega_2} \overline{g} v \text{d}x, \ \text{and} \ \int_{\Gamma} h v_j \text{d}H^{N-1} \rightarrow \int_{\Gamma} h v \text{d}H^{N-1},
\end{equation*}
which combined with \eqref{lim-inf},  \eqref{Fatou1}, \eqref{inf1}, and \eqref{inf2}  provides  $J(v) \leq \liminf\limits_{j \rightarrow +\infty}J(v_j),$ implying that $v$ is a minimizer.
\end{proof}
%\section{The case $h \in L^{\infty}(\Omega) $ }

Now, we are able to regularize  the minimizers of \eqref{func}, given by Theorem \ref{existence-boundness}, under the condition  $h \in L^{\infty}(\Gamma).$

\begin{proof}[Proof of Theorem \ref{unifbound}]
 Let $v \in W^{1,\Phi}_{\varphi}(\Omega)$ be a minimizer of $J$ and $\widehat k_0 \in \mathbb{N}$ such that
\begin{align}\label{defk_0}
\widehat k_0 \geq \| \varphi\|_{L^{\infty}(  \Omega)}.
\end{align} For each $k \geq \widehat k_0$, define the function $v_k : \Omega \rightarrow \mathbb{R} $ given by
\[
v_k =
\begin{cases}
 \mathrm{sgn}(v) k  &\quad\text{if} \ |v|>k,\\
v  &\quad\text{if} \ |v| \leq k,
\end{cases}
\]
where $\mathrm{sgn}(v) = 1$ if $v \geq  0 $, and $\mathrm{sgn}(v)=-1$ if $v<0$. Let $A_k := \{x \in \Omega \ : \ |v(x)| > k\}.$ For each $k \geq \widehat k_0$,
\begin{align*}
	v=v_k \ \text{in} \ A_{k}^c, \quad v_k = \mathrm{sgn}(v) k \  \text{in} \ A_{k}, \ \ \ \text{and} \ \ \ \ (|v|- k)^+\in W^{1,\Phi}_{0}(\Omega). 
\end{align*}
Hence,
\begin{equation}\label{phidifk}
\begin{aligned}
	\int_{\Omega} \Phi (|\nabla v|) -  \Phi(|\nabla v_k|)  dx 
	&= \int_{A_k} \Phi(|\nabla v|) - \Phi(|\nabla v_k|)  dx  + \int_{{A^{c}_{k}}} \Phi(|\nabla v|) - \Phi(|\nabla v_k|) \, \text{d}x
= \int_{A_k} \Phi(|\nabla v|)  \text{d}x.
\end{aligned}
\end{equation}

{We claim that there exists a constant $C>0$ such that
\begin{equation}\label{eq:claim1}
	\int_{A_k} \Phi(|\nabla v| ) \text{d}x \leq C \left( \int_{A_k} (|v|-k)^{1+\alpha} \text{d}x + (k^{1+\alpha}+1) |A_k|\right).
\end{equation}
In order to prove \eqref{eq:claim1} some auxiliary estimates will be needed.}   Since $y $ is a non-negative function, we have that $Y$ is nondecreasing. Thus,
\begin{equation}\label{b1}
	\begin{aligned}
		\int_{\Gamma} Y(v^{+}_{k}) - Y(v^{+}) \text{d}H^{N-1} =& \int_{\Gamma \cap A_k} Y(v^{+}_{k}) - Y(v^{+})  \text{d}H^{N-1} 
		=  \int_{\Gamma \cap A_k \cap \{v\geq 0\}} Y(k) - Y(v) \text{d}H^{N-1}\\  &+ \int_{\Gamma \cap A_k \cap \{v< 0\}} Y((-k)^{+}) - Y(0) \text{d}H^{N-1}   \leq 0.
	\end{aligned}
\end{equation}
Analogously,
\begin{equation}\label{b2}
	\begin{aligned}
		\int_{\Gamma} Y(v^{-}_{k}) - Y(v^{-}) \text{d}H^{N-1} =& \int_{\Gamma \cap A_k} Y(v^{-}_{k}) - Y(v^{-})  \text{d}H^{N-1}  = \int_{\Gamma \cap A_k \cap \{v\geq 0\}} Y(k^{-}) - Y(0)  \text{d}H^{N-1}\\  &+   \int_{\Gamma \cap A_k \cap \{v< 0\}} Y((-k)^{-}) - Y(-v) \text{d}H^{N-1}   %= \int_{\Gamma \cap A_k \cap \{v< 0\}} Y(k) - Y(-v) \text{d}H^{N-1} 
		 \leq 0.
	\end{aligned}
\end{equation}
Moreover,
\begin{equation*}
	\begin{aligned}
		\int_{\Omega_1}  \left( \chi_{\{v_k >0\}} - \chi_{\{v>0\}}\right) \text{d}x =  \int_{\Omega_2}  \left( \chi_{\{v_k >0\}} - \chi_{\{v>0\}}\right) \text{d}x =  \int_{\Gamma}  \left( \chi_{\{v_k >0\}} - \chi_{\{v>0\}}\right) \text{d}H^{N-1} =0,\\
		\int_{\Omega_1}  \left( \chi_{\{v_k \leq0\}} - \chi_{\{v\leq 0\}}\right) \text{d}x =  \int_{\Omega_2}  \left( \chi_{\{v_k \leq 0\}} - \chi_{\{v\leq 0\}}\right) \text{d}x =  \int_{\Gamma}  \left( \chi_{\{v_k \leq 0\}} - \chi_{\{v\leq 0\}}\right) \text{d}H^{N-1} =0,
	\end{aligned}
\end{equation*}
which combined with \eqref{phidifk}, \eqref{b1}, \eqref{b2}, and the inequality $J(v) \leq J(v_k)$ provides that
%\begin{small}
\begin{equation}\label{veryi}
\begin{aligned}
\int_{A_k} \Phi (|\nabla v|) \leq& \int_{\Omega_1} \theta'(F(v^{+}_{k}) - F(v^{+})) \text{d}x+\int_{\Omega_1} \widetilde{\theta}(F(v^{-}_{k}) - F(v^{-})) \text{d}x%+ \int_{\Omega_1}  \left( \chi_{\{v_k >0\}} - \chi_{\{v>0\}}\right) \text{d}x\\+ \int_{\Omega_1} \lambda^{\prime} \left( \chi_{\{v_k \leq 0\}} - \chi_{\{v\leq 0\}}\right) \text{d}x 
+ \int_{\Omega_1} g(v_k - v) \text{d}x\\ &+ \int_{\Omega_2}\sigma' (\overline{F}(v^{+}_{k}) - \overline{F}(v^{+})) \text{d}x+ \int_{\Omega_2} \widetilde{\sigma}(\overline{F}(v^{-}_{k}) - \overline{F}(v^{-})) \text{d}x %+ \int_{\Omega_2}  \left( \chi_{\{v_k >0\}} - \chi_{\{v>0\}}\right) \text{d}x+ \int_{\Omega_2}  \left( \chi_{\{v_k \leq 0\}} - \chi_{\{v\leq 0\}}\right) \text{d}x\\
+\int_{\Omega_2} \overline{g}(v_k - v) \text{d}x\\ &+ \int_{\Gamma}\rho' (Y(v^{+}_{k}) - Y(v^{+})) \text{d}x +\int_{\Gamma} \widetilde{\rho}(Y(v^{-}_{k}) - Y(v^{-})) \text{d}x
%+ \int_{\Gamma} \gamma^{\prime} \left( \chi_{\{v_k >0\}} - \chi_{\{v>0\}}\right) \text{d}x+ \int_{\Gamma} \gamma^{\prime} \left( \chi_{\{v_k \leq 0\}} - \chi_{\{v\leq 0\}}\right) \text{d}x
+\int_{\Gamma}h (v - v_k)\text{d} H^{N-1}\\
\leq& \int_{\Omega_1} \theta'(F(v^{+}_k) - F(v^{+})) \text{d}x +\int_{\Omega_1} \widetilde{\theta}(F(v^{-}_k) - F(v^{-})) \text{d}x+ \int_{\Omega_2} \sigma'(\overline{F}(v^{+}_{k}) - \overline{F}(v^{+}_{k}) )\text{d}x\\&+ \int_{\Omega_2} \widetilde{\sigma}(\overline{F}(v^{-}_{k}) - \overline{F}(v^{-}_{k}) )\text{d}x
%\int_{\Omega_1} g(v_k - v) \text{d}x +\int_{\Omega_2} \overline{g}(v_k - v) \text{d}x +\int_{\Gamma}h (v - v_k)\text{d} H^{N-1} +\gamma\int_{\Gamma} G(v^{+}_{k}) - G(v^{+})\text{d}H^{N-1} + \gamma\int_{\Gamma} G(v^{-}_{k}) - G(v^{-})\text{d}H^{N-1}.
+ \int_{\Omega_1} g(v_k - v) \text{d}x +  \int_{\Omega_2} \bar g(v_k - v) \text{d}x\ +\int_{\Gamma}h (v - v_k)\text{d} H^{N-1}.
\end{aligned}
\end{equation}
%\end{small}
{We now estimate each term of the right hand side of \eqref{veryi}.} Since,
\begin{equation}\label{AB0}
\begin{aligned}
\int_{\Omega_1} |F(v^{+}_k) - F(v^{+})| \text{d}x & = \int_{\Omega_1 \cap A_{k}} |F(v^{+}_k) - F(v^{+}))| \text{d}x = %\int_{\Omega_1 \cap A_{k} \cap \{v> 0\}} |F(v^{+}_k) - F(v^{+}))| \text{d}x = 
\int_{B_k} ( F(|v|) -F(k) ) \text{d}x.
\end{aligned}
\end{equation}
where
%$$A:= \Omega_1 \cap A_k \cap \{v> 0\} \cap \{F(k) > F(|v|)\}, \mbox{ and } 
$B_k:= \Omega_1 \cap A_k \cap \{v> 0\}. $
From \eqref{hipF} and  $(f_1)$ of Lemma \ref{f-properties}, there exist  constants $\alpha, \beta>0$ satisfying
\begin{align*}
	{\max\{1+\delta_0,1+f_0\} <1+\alpha }< (1+\delta_0)^*,
\quad
	{\max\{1+\delta_0,1+\overline{f}_0\} <1+\alpha }< (1+\delta_0)^*, 
\end{align*}
and $ \max\{|\theta_0|, |f_0|, |\bar \theta_0| |\bar f_0|\} < \beta < 1$
such that
\begin{equation*}
	f(t)\leq f(1)\frac{1+f_0}{1+\theta_0}(1+t^\alpha+ t^{-\beta}), \mbox{ and \ }
	\overline{f}(t)\leq \overline{f}(1)\frac{1+\overline{f}_0}{1+\overline{\theta}_0}(1+t^{\alpha}+ t^{-\beta}), \quad \mbox{for all } t>0.
\end{equation*}
From the Mean Value Theorem, there exists $\theta \in (0,1)$, which depends on $x \in \Omega$ and $k \in \mathbb{N}$, such that
\begin{equation}\label{F-1}
\begin{aligned}
	\int_{B_k} (F(|v|) - F(k))\text{d}x&\leq  \int_{B_k} f( (1-\theta)k + \theta |v|)(|v|-k)  \text{d}x\\
	&\leq C \int_{B_k}(1 + |(1-\theta)k + \theta |v||^{-\beta} + |(1-\theta)k + \theta |v||^{\alpha})(|v|-k) \text{d}x \\
	&\leq  C \int_{B_k}(1 + |(1-\theta)k + \theta k|^{-\beta} + |(1-\theta)k + \theta |v||^{\alpha})(|v|-k) \text{d}x \\
	&\leq C \int_{A_k} [1 +  |v|^{\alpha}](|v|-k) \text{d}x ,
\end{aligned}
\end{equation}
where $C>0$ is a constant.
Since $v \in W^{1,\Phi}(\Omega) \hookrightarrow  W^{1,1+\delta_0}(\Omega) $, the embedding $ W^{1,1+\delta_0}(\Omega) \hookrightarrow L^{(1+\delta_0)^*}(\Omega)$ implies that $|v|^{\alpha} \in L^{\frac{1+\alpha}{\alpha}}(\Omega),$ and $(|v|-k) \in L^{1+\alpha}(\Omega).$ Thus, from Young's inequality, we have 
\begin{equation}\label{AB2}
\begin{aligned}
	\int_{A_k}|v|^{\alpha} (|v| -k) \text{d}x &\leq  C \left( \int_{ A_k} (|v|-k + k)^{1+\alpha} \text{d}x + \int_{A_k} (|v| - k)^{1+\alpha} \text{d}x\right)  \\
	&\leq  C \left( \int_{ A_k}  \left( (|v|-k)^{1+\alpha} + k^{1+\alpha}\right) \text{d}x + \int_{ A_k} (|v| - k)^{1+\alpha} \text{d}x\right)  \\
	&\leq  C \left( \int_{ A_k} (|v| - k)^{1+\alpha} \text{d}x + k^{1+\alpha} | A_k| \right).
\end{aligned}
\end{equation}

From \eqref{AB0}, %\eqref{AB1},
 \eqref{F-1}, and \eqref{AB2}  it follows that
\begin{equation}\label{Flast}
 \int_{\Omega_1 } |F(v_k^{+}) - F(v^{+})| \text{d}x
	\leq C\left(  \int_{A_k} (|v|-k) \text{d}x +\int_{ A_k} (|v| - k)^{1+\alpha} \text{d}x + k^{1+\alpha} | A_k| \right).
\end{equation}

Arguing analogously on $v^-$ and $\Omega_2$, we get 
\begin{equation}\label{Flasts}
\begin{aligned}
	 \int_{\Omega_1 } |F(v_k^{-}) - F(v^{-})| \text{d}x
	&\leq C \left( \int_{ A_k} (|v| - k) \text{d}x +\int_{ A_k} (|v| - k)^{1+\alpha} \text{d}x + k^{1+\alpha} |A_k| \right)\\
 \int_{\Omega_2 } |\overline{F}(v_k^{+}) - \overline{F}(v^{+})| \text{d}x
	&\leq C \left( \int_{ A_k} (|v| - k) \text{d}x+\int_{A_k } (|v| - k)^{1+{\alpha}} \text{d}x + k^{1+{\alpha}} |A_k| \right)\\
	\int_{\Omega_2 } |\overline{F}(v_k^{-}) - \overline{F}(v^{-})| \text{d}x&\leq  C  \left( \int_{ A_k} (|v| - k) \text{d}x+\int_{A_k} (|v| - k)^{1+{\alpha}} \text{d}x + k^{1+{\alpha}} |A_k| \right),
\end{aligned}
\end{equation}
for some positive constant $C.$
Since $(|v|-k)^{+} \in W^{1,\Phi}_{0}(\Omega),$  from the continuous embeddings 
$$W^{1,\Phi}_{0}(\Omega) \hookrightarrow  W^{1,1+\delta_0}_{0}(\Omega) \hookrightarrow W^{1,1}_{0}(\Omega) \hookrightarrow L^{\frac{N}{N-1}}(\Omega)\hookrightarrow L^{1}(\Omega),$$ the fact that $\nabla |u| = (\nabla u) \mathrm{sgn}(u)$, {and applying Lemma \ref{tildes} for an arbitrary $\eps>0$,  it follows that
\begin{equation}\label{import}
\begin{aligned}
\int_{A_k} (|v| - k) \text{d}x &= \int_{A_k} (|v|-k)^{+} \text{d}x  \leq C\left( \int_{A_k}  ((|v|-k)^{+})^{\frac{N}{N-1}} \text{d}x\right)^{\frac{N-1}{N}} \leq C\int_{A_k} |\nabla (|v|-k)^{+}| \text{d}x\\
& = C \int_{A_k} |\nabla v| \text{d}x\leq C \left( \varepsilon \int_{A_k} \Phi(|\nabla v|) \text{d}x  + C(\eps)|A_k|\right),
\end{aligned}
\end{equation}}
 with $C$ being a positive constant that does not depend on $\varepsilon,$ and $k.$

Finally, we estimate the last term of the right hand side of \eqref{veryi}, namely
\begin{equation}\label{stop}
\begin{aligned}
\int_{\Gamma} h (v - v_k) \text{d} H^{N-1} &  = \int_{\Gamma \cap A_k \cap \{v> 0\}} h(v-k) \text{d}H^{N-1}  +\int_{\Gamma \cap A_k \cap \{v\leq 0\}} h(v-k) \text{d}H^{N-1} \\
& \leq   2\int_{\Gamma \cap A_k} |h| (|v| - k )\text{d}H^{N-1}.
\end{aligned}
\end{equation}
The inequality \eqref{stop}   combined with  the fact that $(|v|-k)^{+} \in W^{1,\Phi}_{0}(\Omega),$ the embeddings $$W^{1,\Phi}(\Omega_1) \hookrightarrow  W^{1,1+\delta_0}(\Omega_1) \hookrightarrow W^{1,1}(\Omega_1) \hookrightarrow  L^{1}(\Gamma),$$ and \eqref{import} provides 

\begin{equation}\label{imp2}
\begin{aligned}
\int_{\Gamma \cap A_k} |h|(|v| - k) \text{d}H^{N-1} &\leq C \left( \int_{\Gamma \cap A_k} (|v|-k)^{+} \text{d}H^{N-1}\right)\leq C \left( \int_{\Gamma } (|v|-k)^{+} \text{d}H^{N-1}\right)\\ &= C \left( \int_{\Omega_1 \cap A_k} |\nabla (|v| - k )^{+}| \text{d}x + \int_{\Omega_1 \cap A_k} |(|v| - k)^{+}| \text{d}x\right)\\
&\leq C \left( \int_{A_k} |\nabla |v| | \text{d}x + \int_{A_k} (|v| - k) \text{d}x\right) \\
&  \leq C\left( \varepsilon \int_{A_k} \Phi(|\nabla v|) \text{d}x  + C(\eps)|A_k|\right),
\end{aligned}
\end{equation}
where $C>0$  is a constant that does not depend on $\varepsilon$ and $k.$ 

Analogously, we obtain
\begin{equation}\label{diferent}
\begin{aligned}
	\int_{\Omega_1 } g (v_k - v) \text{d}x &=  \int_{ \Omega_1 \cap A_k \cap \{v > 0\}} g (k - v) \text{d}x + \int_{ \Omega_1 \cap A_k \cap \{v\leq0\}} g (-k - v) \text{d}x\leq  2 \int_{\Omega_1 \cap A_k} |g| (|v| -k) \text{d}x.
\end{aligned}
\end{equation}
Using the continuous  embeddings $W^{1,\Phi}(\Omega_1) \hookrightarrow W^{1,1}(\Omega_1)\hookrightarrow L^{\frac{N}{N-1}}(\Omega_1),$ and \eqref{import} we get
\begin{equation*}
\begin{aligned}
	\int_{\Omega_1 \cap A_k} |g| (|v| -k) \text{d}x \leq & C \left(\int_{\Omega_1 \cap A_k} |g|^N \text{d}x\right)^{\frac{1}{N}} \left( \int_{\Omega_1 \cap A_k} (|v| - k)^{\frac{N}{N-1}} \text{d}x\right)^{\frac{N-1}{N}}\\
	\leq & C \left(\int_{\Omega_1} |\nabla (|v|-k)^{+}| \text{d}x +\int_{\Omega_1} |(|v|-k)^{+}| \text{d}x\right) \\
%	=& C_3 \int_{A_k} |\nabla  u| \text{d}x \\
	\leq & C \left( \varepsilon \int_{A_k} \Phi(|\nabla v|) \text{d}x +C(\eps)|A_k|\right),
\end{aligned}
\end{equation*}
where $C>0$ is a constant not depending on $\varepsilon$ and $k,$ which combined with \eqref{diferent} provides that

\begin{equation}\label{dif-3}
	\int_{\Omega_1 } g (v_k - v) \text{d}x \leq 	C \left( {\varepsilon} \int_{A_k} \Phi(|\nabla v|) \text{d}x +C(\eps)|A_k|\right).
\end{equation}
Following the same arguments it follows that
\begin{equation}\label{dif-4}
	\int_{\Omega_2 } \overline{g} (v_k - v) \text{d}x \leq 	C \left( {\varepsilon} \int_{A_k} \Phi(|\nabla v|) \text{d}x +C(\eps)|A_k|\right).
\end{equation}
{Considering  a suitable choice of $\varepsilon$  in \eqref{import}, \eqref{imp2}, \eqref{dif-3}, \eqref{dif-4}, and using \eqref{veryi}, \eqref{Flast}, \eqref{Flasts}, we obtain \eqref{eq:claim1}, where $C= C_\eps$, now depending on the choice of $\eps$, proving the claim. 
Combining \eqref{eq:claim1} with \eqref{auxineq} we have
\begin{equation}\label{fusco}
	\int_{A_k} |\nabla v|^{1+\delta_0} \text{d}x \leq  C_\eps \left(\int_{A_k} (|v|-k)^{1+\alpha} \text{d}x  + (k^{1+\alpha}+1)|A_k| \right).
\end{equation}}

{We claim that  
	\begin{equation}\label{eq:claim2}
	P_{n+1} \leq C D^{n} {P_n}^{1+\xi} \quad \mbox{ for all } n\in\mathbb{N},
		\end{equation} with $C,\xi >0,$ and $D>1 $, being constants which does not depend on $n,$ for sequences 
\begin{align}\label{defK}
P_n:= \int_{A_{k_n}}  (|v|-k_n)^{{1+\alpha}} \text{d}x,  \quad \mbox{and} \quad 	k_n := \frac{K}{2} \left( 1 - \frac{1}{2^{n+1}}\right),
 \ \mbox{ for all } \  n\in\mathbb{N}\cup \{0\},
\end{align}
and for a constant $K>\widehat k_0$, where $\widehat k_0$ is given by \eqref{defk_0}. }

 In fact, using the H\"{o}lder's inequality,  the continuous embedding $W^{1+\delta_0}_0(\Omega) \hookrightarrow L^{(1+\delta_0)^*}(\Omega)$, {and} the fact that $(|u|- k_n)^{+} \in W^{1,1+\delta_0}_0 (\Omega), n\in\mathbb{N},$ we obtain 
\begin{align}\label{fusco2}
	P_{n+1}
	\leq &
	\left(\int_{A_{k_{n+1}}} (|v|- k_{n+1})^{(1+\delta_0)^*} \text{d}x \right)^{{\frac{1+\alpha}{(1+\delta_0)^*}}}{|A_{k_{n+1}}|^{1-\frac{1+\alpha}{(1+\delta_0)^*}}}\notag \\
	\leq & \left(\int_{\Omega}  ((|v|- k_{n+1})^{+})^{(1+\delta_0)^*} \text{d}x\right)^{{\frac{1+\alpha}{(1+\delta_0)^*}}}{|A_{k_{n+1}}|^{1-\frac{1+\alpha}{(1+\delta_0)^*}}} \notag \\
	\leq &  C\left( \int_{\Omega} |\nabla (|v|-k_{n+1})^{+}|^{1+\delta_0} \text{d}x\right)^{ {\frac{1+\alpha}{1+\delta_0}}}{|A_{k_{n+1}}|^{1-\frac{1+\alpha}{(1+\delta_0)^*}}} \notag \\
	=& C\left( \int_{A_{k_{n+1}}} |\nabla v|^{1+\delta_0} \text{d}x\right)^{\frac{{1+\alpha}}{1+\delta_0}}{|A_{k_{n+1}}|^{1-\frac{1+\alpha}{(1+\delta_0)^*}}},
\end{align}
where $C>0$ is a constant that does not depend on $n.$
Combining \eqref{fusco} and \eqref{fusco2}, it follows that
\begin{equation}\label{fusco3}
P_{n+1} \leq C \left(  \int_{A_{k_{n+1}}}(|v|-k_{n+1})^{{1+\alpha}} \text{d}x + (k^{{1+\alpha}}_{n+1} +1)|A_{k_{n+1}}| \right)^{\frac{1+\alpha}{1+\delta_0}}{|A_{k_{n+1}}|^{1-\frac{1+\alpha}{(1+\delta_0)^*}}}.
\end{equation}
Since $k_{n+1} - k_n = \frac{K}{2^{n+3}},$ for all $n\in\mathbb{N}$, we get
\begin{align*}
	\left( \frac{K}{2^{n+3}}\right)^{{1+\alpha}} |A_{k_{n+1}}| = & (k_{n+1} - k_n)^{{1+\alpha}} |A_{k_{n+1}}| \notag 
	= \int_{A_{k_{n+1}}} (k_{n+1} - k_n)^{{1+\alpha}}\text{d}x \notag 
	\leq  \int_{A_{k_{n+1}}} (|v| - k_n)^{{1+\alpha}} \text{d}x,
\end{align*}
which provides
\begin{align}\label{fusco4}
	|A_{k_{n+1}}| \leq  \left(\frac{2^{n+3}}{K}\right)^{{1+\alpha}} 	 \int_{A_{k_{n+1}}} (|v| - k_n)^{{1+\alpha}} \text{d}x 
	\leq   \left(\frac{2^{n+3}}{K}\right)^{{1+\alpha}} P_n.
	%\ {\text{and}~|A_{k_{n+1}}|^{\frac{1+\delta}{ 1+\alpha}-\frac{1+\delta}{(1+\delta_0)^*}} \leq \left(\frac{2^{n+3}}{K}\right)^{{1+\alpha}} J_n}.
\end{align}
Using \eqref{fusco4}, it follows that
\begin{equation}\label{fusco5}
\begin{aligned}
	\int_{A_{k_{n+1}}} (|v|-k_{n+1})^{{1+\alpha}} \text{d}x
&	\leq  {2^\alpha}\int_{A_{k_{n+1}}} (|v|-k_n)^{{1+\alpha}} \text{d}x + {2^\alpha} \int_{A_{k_{n+1}}} (k_{n+1} -k_n)^{{1+\alpha}} \text{d}x  \\
&	\leq {2^\alpha}\int_{A_{k_n}} (|v|-k_n)^{{1+\alpha}} \text{d}x + {2^\alpha} (k_{n+1} -k_n)^{{1+\alpha}} |A_{k_{n+1}}|  \\
&	\leq  {2^\alpha}P_n +{2^\alpha} (k_{n+1} -k_n)^{{1+\alpha}}\left(\frac{2^{n+3}}{K}\right)^{{1+\alpha}} P_n  \\ 
&	=  {2^\alpha}P_n +  {2^\alpha}\left( \frac{K}{2^{n+3}}\right)^{{1+\alpha}} \left(\frac{2^{n+3}}{K}\right)^{{1+\alpha}} P_n   
	=  2^{{1+\alpha}} P_n.
\end{aligned}
\end{equation}

Moreover,
\begin{align}\label{fusco6}
	(k^{{1+\alpha}}_{n+1} +1)|A_{k_{n+1}} | \leq & \left( \left(\frac{K}{2} \left( 1 - \frac{1}{2^{n+2}}\right) \right)^{{1+\alpha}}  +1\right)\left(\frac{2^{n+3}}{K}\right)^{{1+\alpha}}  P_n  \notag \\
	\leq & C 2^{{(1+\alpha)}(n+3)}P_n,
\end{align}
with $C>0$ being a  constant which does not depend on $n.$
Using \eqref{fusco3}, \eqref{fusco5} and \eqref{fusco6}, we obtain   $P_{n+1} \leq C D^n P^{{1+\frac{1+\alpha}{ 1+\delta_0}-\frac{1+\alpha}{(1+\delta_0)^*} }}_n$, where  $C>0$ and $D>1$ are constants does not depending on $n.$ {This proves the claim in \eqref{eq:claim2}, with $\xi:= \frac{1+\alpha}{ 1+\delta_0}-\frac{1+\alpha}{(1+\delta_0)^*}.$}

{Now, we claim that  there exists a constant $K>0$ depending on  $\|v\|_{L^{(1+\delta_0)^*}(\Omega)}$, and consequently,  on $\|v\|_{W^{1,\Phi}(\Omega)}$, such that $|v| \leq \frac{K}{2}$ a.e in $\Omega.$
Indeed, we observe that
	\begin{equation*}
	\begin{aligned}
		%| A_{K_0}|    \left(\frac{K}{4}\right)^{(1+\delta_0)^*} &=
		| A_{k_0}|    k_0  ^{(1+\delta_0)^*} \leq \int_{A_{k_0}}   |v |^{(1+\delta_0)^*}\text{d}x  \leq \int_{\Omega}   |v |^{(1+\delta_0)^*}\text{d}x,
\mbox{	which implies } 
	|A_{k_0}| \leq\left(\frac{1}{k_0}\right)^{(1+\delta_0)^*}  \|v\|_{L^{(1+\delta_0)^*}(\Omega)}^{(1+\delta_0)^*} .
\end{aligned}
\end{equation*}
Hence,
\begin{align*}
	P_0  = &\int_{A_{{k_0}}} ((|v | - k_0)^{+})^{{1+\alpha}}\text{d}x \notag
	\leq \int_{A_{{k_0}}}   |v | ^{{1+\alpha}}\text{d}x \notag
	\leq  \|v\|_{L^{(1+\delta_0)^*}}^{\frac{1+\alpha}{(1+\delta_0)^*}}|A_{k_0}|^{1-\frac{1+\alpha}{(1+\delta_0)^*}}\notag\\
	\leq & \left(\frac{4}{K}\right)^{(1+\delta_0)^*-(1+\alpha)}  \|v\|_{L^{(1+\delta_0)^*}(\Omega)}^{  \frac{1+\alpha}{(1+\delta_0)^*}+   (1+\delta_0)^*-(1+\alpha)   }.
\end{align*}
Thus, it is possible to choose $K$ large enough, depending on  {$\|v\|_{L^{(1+\delta_0)^*}(\Omega)}$},  such that $P_0 \leq C^{-\frac{1}{\xi}} D^{-\frac{1}{\xi^2}}$ . By applying Lemma \ref{lady-aux} it follows that  $P_n \rightarrow 0$. Therefore, we conclude that
\begin{align*}
0 =\lim\limits_{n \rightarrow +\infty} \int_{A_{k_n}} ((|v|-k_n))^{{1+\alpha}} \text{d}x = \int_{A_{\frac{K}{2}}} \left(\left(|v|-\frac{K}{2}\right)\right)^{{1+\alpha}} \text{d}x,
\end{align*}
which provides  $|v| \leq \frac{K}{2}$ a.e in $\Omega,$ proving the claim.}

In what follows we prove the last part of the result, that is, the estimate \eqref{lpart}. In order to obtain such an estimate,    is enough to prove that  $\|v\|_{W^{1,\Phi}(\Omega)}$ depends only on the quantities $\|h\|_{L^{\infty}(\Gamma)}$,  $\theta_0$, $f_0 , \widetilde{\theta}_0 , \widetilde{f}_0, \overline{\theta}_0 , \overline{f}_0$ ,$\delta_0 $, $g_0$, $ \theta',\widetilde{\theta},\sigma',\widetilde{\sigma},\rho',\widetilde{\rho},\widetilde{\lambda},\widetilde{\beta}, \widetilde\gamma, \lambda^{\prime}, \beta^{\prime}, \gamma^{\prime}$, $ \Phi(1)$, $\widetilde{\Phi}(1)$,  $\|g\|_{L^N(\Omega_1)}, \|\overline{g}\|_{L^N(\Omega_2)}$,$\| \varphi\|_{W^{1,\Phi}(\Omega)}$ and the diameter of $\Omega$. In fact, arguing as in the inequality  \eqref{refle1}, there exists $\rho>0$, small enough, and a positive constant $C$ such that, 
$$
\int_{\Omega} \Phi(|\nabla v|) \text{d}x  \leq \rho \int_{\Omega}  \Phi(|\nabla v|) \text{d}x + C\left( \|\varphi \|_{W^{1,1+\delta_0}(\Omega)} + \|\nabla\varphi\|_{L^{1}(\Omega)}+  \|\varphi\|_{L^{1^*}(\Omega)} +1  \right) := \rho \int_{\Omega}  \Phi(|\nabla v|) \text{d}x +A,
$$ which is equivalent to
\begin{align}
	\int_{\Omega} \Phi(|\nabla v|) \text{d}x  \leq&    \frac{A}{1-\rho} := A_\rho. \notag
\end{align}
Since  $v-\varphi\in W^{1,\Phi}_0(\Omega)$, it follows from Lemma \ref{poincare}, and $(\Phi_{4})$ in Lemma \ref{G-properties} that
\begin{equation}\label{2.31}
\begin{aligned}
	\int_{\Omega}  \Phi(|v|)  \text{d}x  \leq&  C\int_{\Omega}  \Phi(|v-\varphi|)  \text{d}x + C\int_{\Omega}  \Phi(|\varphi|)  \text{d}x 
	\leq  C\int_{\Omega}  \Phi(|\nabla v-\nabla \varphi|)  \text{d}x + C\int_{\Omega}  \Phi(|\varphi|)  \text{d}x \\
	\leq& CA_\rho+ C\int_{\Omega}  (\Phi(| \nabla \varphi|)  + \Phi(|\varphi|))  \text{d}x.
\end{aligned}
\end{equation}
Therefore, \[ 	\int_{\Omega} \Phi(|\nabla v|) \text{d}x  +	\int_{\Omega}  \Phi(|v|)  \text{d}x \leq C\left[ \|\varphi \|_{W^{1,1+\delta_0}(\Omega)} + \|\nabla\varphi\|_{L^{1}(\Omega)}+  \|\varphi\|_{L^{1^*}(\Omega)} +1  + \int_{\Omega}  (\Phi(| \nabla \varphi|)  + \Phi(|\varphi|))  \text{d}x\right],\]
where $C>0$ is a constant  depending  only on $\delta_0 , g_0$ and the diameter of $\Omega$. From  \eqref{2.31} and Lemma \ref{modular} we have   that the boundedness of $\|v\|_{W^{1,\Phi}(\Omega)},$ depending only on the desired  quantities mentioned before.
\end{proof}

\begin{remark}
	We point out that the uniform estimate \eqref{lpart}     also hold   true  if the conditions on $F$ and $\overline{F}$ in \eqref{hipF} are  replaced with 
	\begin{equation*}
		\begin{aligned}
			f(t) \leq C(1+ t^{-\beta} + t^{\alpha}),
			\overline{f}(t) \leq C(1+ t^{-\overline{\beta}} + t^{\alpha}),
		\end{aligned}
	\end{equation*}
for all $t>0,$ and
	$$
	{\max\{1+\delta_0,1+f_0, 1+\overline{f}_0\} <1+\alpha }< (1+\delta_0)^*, \beta,\overline{\beta} \in (0,1),
	$$
	with $F^{\prime} (t) = f(t),$ and $\overline{F}^{\prime}(t) = \overline{f}(t), $ $t>0.$

\end{remark}

%%%%%%%%%%%%%%%%%%%%5%%%%%%%%%%%%%%%%%%%%%%%%%%%%%%%%%%%%%%%%%%%%%%%%%%%%%%%%%%%%%%%%%%%%%%%%%%%%%%%%%%%%%%%%%%%%%%%%%%%%%%%%%

\section{Morrey  type results for  the gradient   and H\"{o}lder estimates}
\hfill

The purpose of this section is to obtain    Morrey  type regularity results for the gradient of the  minimizers of  \eqref{func}.    
In what follows, we present the proof of Theorem \ref{morreyg},  based on the proof of \cite[Proposition 4.3]{queirozt}.

\begin{proof}[Proof of Theorem \ref{morreyg}]
Consider an open set $\Omega^{\prime} \subset \subset \Omega , x_0 \in \Omega^{\prime},$ and open balls  $B_r (x_0) \subset B_R (x_0) \subset \subset \Omega$. Let  $v^* \in W^{1,\Phi}(B_R(x_0))$ be the $\Phi-$harmonic replacement of the minimizer $v$ in $B_R(x_0),$ that is, $v^*$ satisfies \eqref{hreplace}. 
From the Maximum Principle \cite[Lemma 2.8]{MW}  we have $\| v^*\|_{L^{\infty} (B_R)} \leq \| v\|_{L^{\infty}(\Omega)}.$ {Thus,  from Lemma \ref{iotal0} we have that
\begin{equation}\label{hreplacement}
\int_{B_r} \Phi (|\nabla v^*|)\text{d}x \leq C \left( \frac{r}{R}\right)^N \int_{B_R }\Phi (|\nabla v^*|)\text{d}x, \ \mbox{for} \ 0<r \leq R,
\end{equation}
where $C=C(N,\delta_0,g_0, \| v\|_{L^{\infty}(\Omega)})$. In addition,  from $(\Phi_4)$ in Lemma \ref{G-properties} we obtain}

\begin{equation}\label{m1.0}
\begin{aligned}
\int_{B_r} \Phi(|\nabla v|) \text{d}x &\leq \int_{B_r}\Phi (|\nabla v - \nabla v^*| + |\nabla v^*|) \text{d}x \\
 & \leq C \int_{B_r} \Phi (|\nabla v - \nabla v^*|) \text{d}x + C \int_{B_r} \Phi (|\nabla v^*|) \text{d}x \\
 & \leq C  \int_{B_r} \Phi (|\nabla v - \nabla v^*|) \text{d}x + C  \left( \frac{r}{R}\right)^N \int_{B_R} \Phi (|\nabla v^*|)\text{d}x \\
 & \leq C  \int_{B_r} \Phi (|\nabla v - \nabla v^*|) \text{d}x + C  \left( \frac{r}{R}\right)^N \int_{B_R} \Phi (|\nabla v|)\text{d}x,
\end{aligned}
\end{equation}
where in the last inequality we applied the minimizing property of the Dirichlet integral  $\int_{B_R} \Phi (|\nabla v^*|)\text{d}x $.
From Lemma \ref{liota}, and the inequality $ab \leq \frac{1}{2}(a^2 + b^2), $ it follows that for each $\mu \in(0,N),$ there is a constant $C$, depending on $\mu$, such that
\begin{equation}\label{m2.0}
\begin{aligned}
	\int_{B_{R}}\Phi(|\nabla v-\nabla v^*|)\text{d}x&\leq  C\int_{B_{R}}(\Phi(|\nabla v|)-\Phi(|\nabla v^*|))\text{d}x +CR^{\frac{\mu}{2}}\bigg(\int_{B_{R}}(\Phi(|\nabla v|)-\Phi(|\nabla v^*|))\text{d}x\bigg)^{\frac{1}{2}}\\
&\leq 
 C\int_{B_{R}}(\Phi(|\nabla v|)-\Phi(|\nabla v^*|))\text{d}x\ + CR^{\mu} ,
\end{aligned}
\end{equation}
where $C>0$ is a   constant, which does not depend on $v,r,$ $R$.
Consider that $v^* \equiv v$ in ${B_R(x_0)}^c \cap \Omega.$ Since, $v$ is a minimizer we obtain that
\begin{equation}\label{m3.0}
\begin{aligned}
&\int_{B_{R}}(\Phi(|\nabla v|)-\Phi(|\nabla v^*|))\text{d}x\leq \int_{\Omega_1 \cap B_R} W(v^*) - W(v) + g (v^* - v) \text{d}x + \int_{\Omega_2 \cap B_R} \overline{W}(v^*) - \overline{W}(v) + \overline{g} (v^* - v) \text{d}x \\
&+ \int_{\Gamma \cap B_R} h (v - v^*) \text{d}H^{N-1} + \int_{\Gamma \cap B_R}   \rho'(Y(({v^*)}^{+}) - Y(v^{+})) \text{d}H^{N-1} + \int_{\Gamma \cap B_R }\widetilde{\rho}(Y(({v^*)}^{-}) - Y(v^{-})) \text{d}H^{N-1}\\
& +\int_{\Gamma \cap B_R}  \gamma^{\prime} \left( \chi_{\{v^*>0\}} -\chi_{\{v>0\}}\right) \text{d}H^{N-1} +\int_{\Gamma \cap B_R}  \widetilde\gamma \left( \chi_{\{v^*\leq 0\}} -\chi_{\{v\leq0\}}\right) \text{d}H^{N-1}.
\end{aligned}
\end{equation}
Thus, from \eqref{m1.0}-\eqref{m3.0} we obtain that

\begin{equation}\label{mkey}
\begin{aligned}
\int_{B_r}	\Phi (|\nabla v|) \text{d}x \leq & \, C R^{\mu} + C\left( \frac{r}{R}\right)^N \int_{B_R} \Phi (|\nabla v|)\text{d}x +  \int_{\Gamma \cap B_R } \rho' ( Y((v^*)^{+}- Y(v^{+})) \text{d}H^{N-1} \\
&+  \int_{\Gamma \cap B_R} \widetilde{\rho}(Y({(v^*)}^{-}) - Y(v^{-})) \text{d}H^{N-1} +\int_{\Gamma \cap B_R }  \gamma^{\prime} \left( \chi_{\{v^*>0\}} -\chi_{\{v>0\}}\right) \text{d}H^{N-1} \\
& +\int_{\Gamma \cap B_R} \widetilde{\gamma} \left( \chi_{\{v^*\leq 0\}} -\chi_{\{v\leq0\}}\right) \text{d}H^{N-1} + \int_{\Gamma \cap B_R} h (v - v^*) \text{d}H^{N-1}\\
&+\int_{\Omega_1 \cap B_R} g(v^* - v) \text{d}x +  \int_{\Omega_2 \cap B_R} \overline{g}(v^* - v) \text{d}x \\
&  + \int_{\Omega_1 \cap B_R} W(v^*) - W(v) \text{d}x  + \int_{\Omega_2 \cap B_R} \overline{W}(v^*) - \overline{W}(v) \text{d}x .
\end{aligned}
\end{equation}
Note that
\begin{equation}\label{boundgprime}
	\int_{\Gamma \cap B_R } \gamma^{\prime} \left( \chi_{\{v^* >0\}} - \chi_{\{v >0\}} \right)\text{d}H^{N-1} +\int_{\Gamma \cap B_R} \widetilde\gamma \left( \chi_{\{v^*\leq 0\}} - \chi_{\{v \leq 0\}} \right)\text{d}H^{N-1}  \leq CR^{N-1},
\end{equation}
for some constant $C>0$, depending on $\gamma'$ and $\widetilde\gamma.$ Moreover, the continuous embedding $W^{1,1}(\Omega_1) \hookrightarrow L^{1}(\Gamma),$ the fact that $\| v^*\|_{L^{\infty} (B_R)} \leq \| v\|_{L^{\infty}(\Omega)}$, and Lemma \ref{tildes}  provides
\begin{equation}\label{bgood}
\begin{aligned}
\int_{\Gamma \cap B_R} h (v - v^*) \text{d}H^{N-1}  &\leq \| h\|_{L^{\infty}(\Gamma)} \int_{\Gamma \cap B_R} |v- v^*| \text{d}H^{N-1} \leq \| h\|_{L^{\infty}(\Gamma)} \int_{\Gamma } |v- v^*| \text{d}H^{N-1}\\
& \leq \| h\|_{L^{\infty}(\Gamma)} \left( \int_{\Omega_1} |v -v^*| \text{d}x + \int_{\Omega_1} |\nabla v -\nabla v^*| \text{d}x  \right)\\
& = \| h\|_{L^{\infty}(\Gamma)} \left( \int_{\Omega_1 \cap B_R} 2\|v \|_{L^\infty(\Omega)} \text{d}x + \int_{\Omega_1 \cap B_R} |\nabla v -\nabla v^*| \text{d}x  \right) \\
%& \leq  \| h\|_{L^{\infty}(\Gamma)} \left( \int_{ B_R} |v -v^*| \text{d}x + \int_{B_R} |\nabla v -\nabla v^*| \text{d}x  \right)\\
& \leq CR^N + \varepsilon \int_{B_R} \Phi (|\nabla v - \nabla v^* |) \text{d}x,
\end{aligned}
\end{equation}
where $\varepsilon > 0$ and $C>0$ is  a constant.  Furthermore, since  $v^* \equiv v$ in ${B_R(x_0)}^c \cap \Omega$, it follows that $v^* -v \in W^{1,1}_{0}(\Omega)$, hence, in view of the Sobolev embedding $W^{1,1}_{0}(\Omega) \hookrightarrow L^{\frac{N}{N-1}}(\Omega),$ Lemma \ref{tildes}, and the fact that $g \in L^{q}(\Omega_1)$, and $\overline{g} \in L^{\overline{q}}(\Omega_2), $ with $q,\overline{q}\geq N$, we obtain
\begin{equation}\label{overb}
	\begin{aligned}
		\int_{\Omega_1 \cap B_R } g (v^* -v)\text{d}x +  \int_{\Omega_2 \cap B_R } \overline{g} (v^* -v)\text{d}x  \leq & \left[\left( \int_{B_R} |g|^{N}\text{d}x\right)^{\frac{1}{N}} + \left( \int_{B_R} |\overline{g}|^{N}\text{d}x\right)^{\frac{1}{N}}\right]\left( \int_{B_R} |v^* -v|^{\frac{N}{N-1}} \text{d}x\right)^{\frac{N-1}{N}} \\
		%&  \left( \int_{B_R} |v^* -v|^{\frac{N}{N-1}} \text{d}x\right)^{\frac{N-1}{N}} \\
		\leq\, & C \int_{B_R} |\nabla v^* - \nabla v |\text{d}x
		\leq  C\left( R^{N} + \varepsilon \int_{B_R} \Phi(|\nabla v^* - \nabla v|)\text{d}x \right),
	\end{aligned}
\end{equation} 
for an arbitrary $\eps >0$, some constant $C>0,$ and  $R>0.$ Additionally, 
in view  of Lemma \ref{f-properties}, item $(F_1)$, and using that $v, v^* \in L^\infty(\Omega),$ it follows that
\begin{equation}\label{wb}
\int_{\Omega_1 \cap B_R} W(v^*) - W(v) \text{d}x + \int_{\Omega_2 \cap B_R} \overline{W}(v^*) - \overline{W}(v) \text{d}x \leq C R^N,
\end{equation}
where we applied the same argument as in \eqref{boundgprime} to control the terms involving the characteristic functions, hence, $C>0$ is a constant depending on $\lambda',\widetilde\lambda, \beta', \widetilde\beta$. 

From now on the proof will be {split} in some cases.

\textbf{Case 1: $y=Y^{\prime}$  is bounded.}
Define the set $A:= \Gamma \cap B_R (x_0) \cap (\{(v^*)^{+}>0\} \cup \{v^{+}>0\})$. From the Mean Value Theorem, we have that there exists $\theta \in (0,1)$, which depends on $x \in \Omega$ and $R$, such that
\begin{equation}\label{meanb}
	\begin{aligned}
		\int_{\Gamma \cap B_R} Y((v^*)^{+}) - Y(v^{+}) \text{d}H^{N-1} &=	\int_{A} Y((v^*)^{+}) - Y(v^{+}) \text{d}H^{N-1} \\
		&\leq  	\int_{A}  |y((1-\theta)v^{+} + \theta (v^*)^{+})| | (v^*)^{+}- v^{+}| \text{d}H^{N-1}.\\
	\end{aligned}
\end{equation}
From Lemma \ref{ineqpart}, which also holds in $\Gamma \cap B_R (x_0)$, we have
\begin{equation}\label{ol}
\int_{\Gamma \cap B_R} |(v^*)^{+} - v^{+}|\text{d}H^{N-1} \leq \int_{\Gamma \cap B_R} |v^* - v| \text{d}H^{N-1}.
\end{equation}
Since $y$ is bounded, from  the reasoning of  \eqref{bgood}, \eqref{meanb}, and \eqref{ol}, it follows that

\begin{equation}\label{boundary1}
	\begin{aligned}
		\int_{\Gamma \cap B_R} Y((v^*)^{+}) - Y(v^{+}) \text{d}H^{N-1} &\leq C \int_{\Gamma \cap B_R} | v^* - v| \text{d}H^{N-1} \leq  CR^N + \varepsilon \int_{B_R} \Phi (|\nabla v - \nabla v^* |) \text{d}x,
	\end{aligned}
\end{equation}
where $C>0$ is a constant, and $\varepsilon>0 $ is an arbitrary constant. Similarly, we have

\begin{equation}\label{boundary2}
\begin{aligned}
	\int_{\Gamma \cap B_R} Y((v^*)^{-}) - Y(v^{-}) \text{d}H^{N-1}  &\leq  C \int_{\Gamma \cap B_R} | v^* - v| \text{d}H^{N-1} \leq CR^N + {\varepsilon} \int_{B_R} \Phi (|\nabla v - \nabla v^* |) \text{d}x,
\end{aligned}
\end{equation}
where $C>0$ is a constant, and $\varepsilon>0$ is an arbitrary constant.
From \eqref{mkey}-\eqref{boundgprime}, and a suitable choice of $\varepsilon>0$ in \eqref{bgood}-\eqref{overb}, and \eqref{boundary1}-\eqref{boundary2} respectively, we have 
\begin{equation}\label{gen}
\begin{aligned}
\int_{B_r} \Phi(|\nabla v|)\text{d}x & \leq C\left( \frac{r}{R}\right)^{N}\int_{B_R} \Phi (|\nabla v|)\text{d}x  +  CR^{\mu}  + (\gamma'+\widetilde\gamma )CR^{N-1} +  CR^{N} \\
&\leq C\left( \frac{r}{R}\right)^{N}\int_{B_R} \Phi (|\nabla v|)\text{d}x  +  CR^{\mu}  +(\gamma'+\widetilde\gamma ) CR^{N-1},
\end{aligned}
\end{equation}
for some constants $C>0,$ which does not depend on $R>0.$ 

\textbf{Case 2: $y=Y^{\prime}$ is nonincreasing.} \
Define the sets $ A_1 :=  \Gamma \cap B_R (x_0) \cap \{ (v^*)^{+} \geq  v^{+}  \},$ and
 ${A_2 :=  \Gamma \cap B_R (x_0) \cap \{ (v^*)^{+} < v^{+}  \}.}$
From Remark \ref{similar} we have 
\begin{equation}\label{case21}
\begin{aligned}
	\int_{\Gamma \cap B_R} Y((v^*)^{+}) - Y(v^{+}) \text{d}H^{N-1} & = \int_{A_1} Y((v^*)^{+}) - Y(v^{+}) \text{d}H^{N-1}+ \int_{A_2} Y((v^*)^{+}) - Y(v^{+}) \text{d}H^{N-1}\\
	&\leq \int_{A_1} Y((v^*)^{+} - v^{+}) \text{d}H^{N-1}  + \int_{A_2} Y(v^{+} - (v^*)^{+})\text{d}H^{N-1}\\
	& \leq Y(1) \max\left\{ \int_{\Gamma \cap B_R} |(v^*)^{+} - v^{+}|^{1+\widetilde{\theta}_0}\text{d}H^{N-1}, \int_{\Gamma \cap B_R} |(v^*)^{+} - v^{+}|^{1+\widetilde{f}_0}\text{d}H^{N-1}\right\}.
\end{aligned}
\end{equation}
If $-1 < \widetilde\theta_0 \leq \widetilde f_0 < 0$, we obtain $|t|^{1 + \widetilde{\theta}_0} \leq 1 + |t|$ and $|t|^{1 + \widetilde{f}_0} \leq 1 + |t|$ for all $t \in \mathbb{R}$, which, combined with the reasoning of \eqref{bgood}, \eqref{ol},  and \eqref{case21} provides that
\begin{equation}\label{nafil1}
\begin{aligned}
	\int_{\Gamma \cap B_R} Y((v^*)^{+}) - Y(v^{+}) \text{d}H^{N-1} & \leq C \int_{\Gamma \cap B_R } 1+ |(v^*)^{+} - v^{+}| \text{d}H^{N-1}  \leq C \int_{\Gamma \cap B_R } 1+ |v^*- v| \text{d}H^{N-1}\\
	& \leq C\left(R^{N-1} + R^{N} +\varepsilon\int_{B_R} \Phi (|\nabla v - \nabla v^* |) \text{d}x\right),
\end{aligned}
\end{equation}
where $C>0$ is a constant, and $\varepsilon>0$ is an arbitrary constant. A similar argument
gives 
\begin{equation}\label{nafil2}
\begin{aligned}
\int_{\Gamma \cap B_R} Y((v^*)^{-}) - Y(v^{-}) \text{d}H^{N-1}
 &\leq C\left(R^{N-1} + R^{N} + {\varepsilon}\int_{B_R} \Phi (|\nabla v - \nabla v^* |) \text{d}x\right),
\end{aligned}
\end{equation}
where $C>0$ is a constant, and ${\varepsilon}>0$ is an arbitrary constant.
From \eqref{mkey}-\eqref{boundgprime}, and a suitable choice of $\varepsilon>0$ in \eqref{bgood}-\eqref{overb}, and \eqref{nafil1}-\eqref{nafil2}, respectively, provides
\begin{equation*}
	\int_{B_r} \Phi(|\nabla v|)\text{d}x \leq C\left( \frac{r}{R}\right)^{N}\int_{B_R} \Phi (|\nabla v|)\text{d}x  +  CR^{\mu} + (\gamma+ \gamma^{\prime}) CR^{N-1} +  CR^{N}.
\end{equation*}
Suppose that  $-1<\widetilde{\theta}_0 <0\leq \widetilde{f}_0.$ Since $1\leq 1+\widetilde{f}_0 \leq (1+\delta_0)_{*},$ and $\| v^*\|_{L^{\infty}(B_R )} \leq \| v\|_{L^{\infty}(\Omega)}$ we obtain that $|v^* - v|^{1 + \widetilde{f}_0} \in W^{1,1}(\Omega_1).$ Thus, from the Trace Theorem and Lemma \ref{tildes}, we obtain that
{\small\begin{equation}\label{bgood2}
\begin{aligned} 
\int_{\Gamma \cap B_R } |(v^*)^{+} - v^{+}|^{1+\widetilde{f}_0}\text{d}H^{N-1} &  \leq \int_{\Gamma} |v^* - v|^{1+\widetilde{f}_0}\text{d}H^{N-1} \leq C\left( \int_{B_R} |v^* - v|^{1+\widetilde{f}_0} \text{d}x + \int_{B_R} (1+\widetilde{f}_0) |v^* - v|^{\widetilde{f}_0} |\nabla v^* - \nabla v| \text{d}x\right)\\
& \leq C\left( R^{N} + \int_{B_R} |\nabla v^* - \nabla v| \text{d}x \right) \leq C \left(R^{N} + \varepsilon\int_{B_R} \Phi (|\nabla v^* - \nabla v|)\text{d}x\right),
\end{aligned}
\end{equation}}
where $C>0$ is a constant, and  $\varepsilon >0 $ is an arbitrary constant. Then, \eqref{case21}, \eqref{bgood2}, and the inequality $|t|^{1 + \widetilde{\theta}_0} \leq 1 + |t|,$ for all $t \in \mathbb{R}$ provides 
\begin{equation}\label{case22}
\begin{aligned}
	\int_{\Gamma \cap B_R} Y((v^*)^{+}) - Y(v^{+}) \text{d}H^{N-1} \leq C \left( R^{N-1}+R^{N} + \varepsilon\int_{B_R  } \Phi (|\nabla v^*- \nabla v |)\text{d}x\right),
\end{aligned}
\end{equation}
where $C>0$ is a constant, and  $\varepsilon>0 $ is an arbitrary constant. Similarly,

\begin{equation}\label{bgood3}
	\begin{aligned}
		\int_{\Gamma \cap B_R} Y((v^*)^{-}) - Y(v^{-}) \text{d}H^{N-1} \leq C \left( R^{N-1}+R^{N} + {\varepsilon}\int_{B_R} \Phi (|\nabla v^*- \nabla v |)\text{d}x\right),
	\end{aligned}
\end{equation}
where $C>0$ is a constant, and  ${\varepsilon} >0 $ is an arbitrary constant.
Thus, from  \eqref{mkey}-\eqref{boundgprime} and a suitable choice of $\varepsilon>0$ in \eqref{bgood}-\eqref{overb}, and \eqref{bgood2}-\eqref{bgood3}, respectively, we obtain that 

\begin{equation*}
	\int_{B_r} \Phi(|\nabla v|)\text{d}x \leq C\left( \frac{r}{R}\right)^{N}\int_{B_R} \Phi (|\nabla v|)\text{d}x  +  CR^{\mu} + (\widetilde\gamma+ \gamma^{\prime}) CR^{N-1} +  CR^{N},
\end{equation*}
for some constant $C>0,$ which does not depend on $R>0,$ and $\mu \in (0,N)$, which was fixed previously.
Now, suppose that  $0\leq \widetilde{\theta}_0\leq \widetilde{f}_0.$ Using \eqref{case21} and arguing as in \eqref{bgood2}, it follows that
\begin{equation}\label{cach1}
\begin{aligned}
\int_{\Gamma \cap B_R} Y((v^*)^{+}) - Y(v^{+}) \text{d}H^{N-1} 
\leq\, & C\left( \int_{B_R} |v^* - v|^{1+\widetilde{\theta}_0} \text{d}x + \int_{B_R} (1+\widetilde{\theta}_0) |v^* - v|^{\widetilde{\theta}_0} |\nabla v^* - \nabla v| \text{d}x\right) \\
&+ C\left( \int_{B_R} |v^* - v|^{1+\widetilde{f}_0} \text{d}x + \int_{B_R} (1+\widetilde{f}_0) |v^* - v|^{\widetilde{f}_0} |\nabla v^* - \nabla v| \text{d}x\right)\\
\leq\, & C\left( R^{N} + \varepsilon\int_{B_R} \Phi(|\nabla v - \nabla v^*|) \text{d}x \right),
\end{aligned}
\end{equation}
where $C>0$ is a constant, and $\varepsilon>0$ is an arbitrary constant. Similarly,

\begin{equation}\label{cach2}
\begin{aligned}
\int_{\Gamma \cap B_R} Y((v^*)^{-}) - Y(v^{-}) \text{d}H^{N-1}
&\leq C \left( R^{N} + {\varepsilon}\int_{B_R } \Phi (|\nabla v^*- \nabla v |)\text{d}x\right),
\end{aligned}
\end{equation}
where $C>0$ is a constant, and ${\varepsilon}>0$ is an arbitrary constant. 
Thus, it follows  from  \eqref{mkey}-\eqref{boundgprime} and a suitable choice of $\varepsilon>0$ in \eqref{bgood}-\eqref{overb}, and \eqref{cach1}-\eqref{cach2}, respectively, that 

\begin{equation*}
	\int_{B_r} \Phi(|\nabla v|)\text{d}x \leq C\left( \frac{r}{R}\right)^{N}\int_{B_R } \Phi (|\nabla v|)\text{d}x  +  CR^{\mu} +  (\gamma'+\widetilde\gamma) CR^{N-1} +  CR^{N},
\end{equation*}
for some constant $C>0,$ which does not depend on $R>0,$ and $\mu \in (0,N)$, which was fixed before.

\textbf{Case 3: $y=Y^{\prime}$ is nondecreasing.}
Consider the sets $A_1$,  $A_2$ as before, and define $$\widetilde A:= A_1 \cap (\{(v^*)^{+} >0\} \cup \{v^{+} > 0\}).$$
 We have $\| (v^*)^{+}\|_{L^{\infty}(B_R(x_0))} \leq \| v\|_{L^{\infty}(\Omega)}. $ From Remark \ref{similar}, it follows that
\begin{equation}\label{c1}
\begin{aligned}
\int_{\Gamma \cap B_R} Y((v^*)^{+}) - Y(v^{+}) \text{d}H^{N-1} &= \int_{A_1} Y((v^*)^{+}) - Y(v^{+})  \text{d}H^{N-1} +\int_{A_2} Y((v^*)^{+}) - Y(v^{+})  \text{d}H^{N-1}\\
& \leq %\int_{A_1} Y((v^*)^{+}) - Y(v^{+})  \text{d}H^{N-1} =
 \int_{\widetilde A} Y((v^*)^{+}) - Y(v^{+})  \text{d}H^{N-1}
 \leq  \int_{\widetilde A} y(2\|v \|_{L^\infty(\Omega)})((v^*)^{+} - v^{+})\text{d}H^{N-1}\\
& \leq C \int_{\Gamma \cap B_R} |(v^*)^{+} - v^{+}| \text{d}H^{N-1}\leq C \int_{\Gamma \cap B_R} |v^* - v|\text{d}H^{N-1}\\
& \leq C \left( R^{N} + {\varepsilon}\int_{\widetilde A} \Phi(|\nabla v - \nabla v^*|)\right)\text{d}x,
\end{aligned}
\end{equation}
for some constant $C>0,$ and ${\varepsilon} >0$ arbitrary, where it was used the reasoning of \eqref{bgood}, and Lemma \ref{tildes}. Analogously, 
\begin{equation}\label{c2}
\begin{aligned}
\int_{\Gamma \cap B_R} Y((v^*)^{-}) - Y(v^{-}) \text{d}H^{N-1} & \leq  C \int_{\Gamma \cap B_R} |v^* - v|\text{d}H^{N-1}\leq C \left( R^{N} + \varepsilon\int_{B_R} \Phi(|\nabla v - \nabla v^*|)\text{d}x\right),
\end{aligned}
\end{equation}
for some constant $C>0,$ and $\varepsilon>0$ arbitrary.
From \eqref{mkey}-\eqref{boundgprime}, and a suitable choice of $\varepsilon>0$ in \eqref{bgood}-\eqref{overb}, and \eqref{c1}-\eqref{c2}, respectively, we obtain  

\begin{equation*}
	\int_{B_r} \Phi(|\nabla v|)\text{d}x \leq C\left( \frac{r}{R}\right)^{N}\int_{B_R} \Phi (|\nabla v|)\text{d}x  +  CR^{\mu} + (\gamma^{\prime}+\widetilde\gamma) CR^{N-1} +  CR^{N},
\end{equation*}
for some constant $C>0$, which does not depend on $R>0,$ and $\mu \in (0,N)$, which was chosen previously.

Now, we conclude the proof. In order to prove the item  $(a)$,  we consider the Case 1. Recall that \eqref{gen} was obtained for a fixed $\mu \in (0,N)$. Thus, \eqref{gen} also holds for  $\mu = N-1.$ Furthermore, we consider $R_0 >0$ such that $0<R_0 < \text{dist}(\partial \Omega, \partial \Omega^{\prime}) /2,$ where $\Omega^{\prime} \subset \subset \Omega ,\ x_0 \in \Omega^{\prime},$ as mentioned before. Since $\mu = N-1,$ and  $\gamma^{\prime} > 0$ or $\widetilde \gamma>0$, from \eqref{gen}  we get
\begin{equation}\label{m1}
	\int_{B_r}  \Phi(|\nabla v|)\text{d}x \leq C\left( \frac{r}{R}\right)^{N}\int_{B_R}  \Phi (|\nabla v|)\text{d}x  +  CR^{N-1},
\end{equation}
with $C>0$ a constant, for all $0<r < R \leq R_0 .$
Fix $\rho \in (0,1),$ and consider $\tau \in (0,1) $ such that $\tau^{-\rho} > C,$ where $C$ is the constant of \eqref{m1}. Define $r:= \tau R.$ From \eqref{m1} we have
\begin{equation}\label{m2}
\begin{aligned}
	\int_{B_{\tau R}} \Phi(|\nabla v|)\text{d}x&	\leq C \tau^{N}\int_{B_R} \Phi (|\nabla v|)\text{d}x  +CR^{N-1}  \leq \tau^{N-\rho} \int_{B_R} \Phi (|\nabla v|)\text{d}x   + C R^{N-1},
\end{aligned}
\end{equation}
for each $0 < R < R_0$. Thus, from  Lemma \ref{glemma} it follows that
\begin{equation}\label{mf}
\int_{B_r} \Phi(|\nabla v|) \text{d}x \leq C \left( \left( \frac{r}{R_0}\right)^{N-1} \int_{B_{R_0}} \Phi (|\nabla v|)\text{d}x + r^{N-1}\right),
\end{equation}
for all $0<r<R\leq R_0,$ where $C$ is a positive constant that does not depend on $r$ and $R.$ % depends only on $n, \|v \|_{L^{\infty}(\Omega)}, \tau , \rho, \delta_0, g_0, \theta_0 , f_0 , \widetilde{\theta_0}, \widetilde{f_0}, \overline{\theta}_0,$ and $\overline{f}_0.$  
Now, we consider $\Omega^{\prime \prime}$ an open domain such that
$$
\Omega^{\prime} \subset \bigcup_{x_0 \in \Omega^{\prime} } B_{R_0}(x_0) \subset \subset \Omega'' \subset \subset \Omega.
$$
For $0<r<R_0,$ and $x_0 \in \Omega^{\prime}$,  from \eqref{mf} it follows that

\begin{equation}\label{mg} \int_{B_r \cap \Omega^{\prime}} \Phi(|\nabla v|)\text{d}x \leq C \left( \left( \frac{r}{R_0}\right)^{N-1} \int_{\Omega^{\prime \prime}} \Phi(|\nabla v|)\text{d}x + r^{N-1}\right).\end{equation}
Thus, for an arbitrary $0 \leq  \zeta \leq  N-1 $, it yields
\begin{equation*}
\begin{aligned}
r^{-\zeta} \int_{B_r  \cap \Omega^{\prime}} \Phi(|\nabla v|) \text{d}x &\leq C \left( \frac{r^{N-1 -\zeta}}{{R_{0}}^{N-1}} \int_{\Omega^{\prime \prime}} \Phi (|\nabla v|) \text{d}x+ r^{N-1 - \zeta}\right)  \leq \text{diam}(\Omega^{\prime \prime})^{N-1-\zeta}\left( \frac{1}{{R_0}^{N-1}} \int_{\Omega^{\prime \prime}} \Phi(|\nabla v|) \text{d}x + 1 \right),
\end{aligned}
\end{equation*}
for all $0<r<R_0,$ finishing the proof of item $(a)$. The items $(c)$ and $(e)$ follow the same arguments assuming the Cases 2 and 3, respectively. In order to prove the item  $(b),$  we consider the Case 1 again, and the fact that \eqref{gen} holds for each  $\mu \in (0,N).$ Moreover, in this item, the sum $(\gamma' + \widetilde\gamma)$ is null, hence, we are able to repeat the arguments in \eqref{m1}-\eqref{mg} replacing $N-1$ by $\mu$. Therefore, for an arbitrary $0 \leq  \zeta \leq  \mu $, it yields
\begin{equation*}
	\begin{aligned}
		r^{-\zeta} \int_{B_r  \cap \Omega^{\prime}} \Phi(|\nabla v|) \text{d}x &\leq C \left( \frac{r^{\mu -\zeta}}{{R_{0}}^{\mu}} \int_{\Omega^{\prime \prime}} \Phi (|\nabla v|) \text{d}x+ r^{\mu - \zeta}\right)  \leq \text{diam}(\Omega^{\prime \prime})^{\mu-\zeta}\left( \frac{1}{{R_0}^{N-1}} \int_{\Omega^{\prime \prime}} \Phi(|\nabla v|) \text{d}x + 1 \right),
	\end{aligned}
\end{equation*}
for all $0<r<R_0,$ finishing the proof of item $(b)$. The items $(d)$ and $(f)$ follow the same arguments considering the Cases 2 and 3, respectively, which completes the proof.
\end{proof}

Now, we are able to present the proof of Corollary \ref{holderg}.

\begin{proof}[Proof of Corollary \ref{holderg}]
	Consider $v \in W^{1,\Phi}_{\varphi}(\Omega)$  a minimizer of \eqref{func}. In order to prove the item $(a)$, note that from Theorem \ref{embedding} it follows that $ v \in W^{1,1+\delta_0}(\Omega).$  Fix an open domain $\Omega^{\prime} \subset \subset \Omega$ with $\text{diam} (\Omega^{\prime})\leq 1$, and an open ball $B_R(x_0) \subset \Omega^{\prime}.$ Let $\zeta:= N- (1+\delta_0) + \sigma(1+\delta_0),$ where $0<\sigma\leq  \frac{\delta_0}{1 + \delta_0}.$ Note that $0<\zeta \leq N-1. $ Define
	$v_{x_0 , R} := \fint_{B_R } v \text{d}x.$
Applying  Poincaré's inequality \cite[Theorem 1.51]{maly}, using \eqref{auxineq}, and the item $(a)$ of Theorem \ref{morreyg}, it follows that
	
	\begin{equation*}
		\begin{aligned}
			\int_{B_r} |v - v_{x_0 , r}|^{1+\delta_0} \text{d}x &\leq C R^{1+\delta_0} \int_{B_R} |\nabla v|^{1+\delta_0}\leq C R^{1+\delta_0} \left( \int_{B_R} \Phi (|\nabla v|) \text{d}x + \int_{B_R} 1 \text{d}x\right) \\
			&\leq C R^{1+\delta_0} (R^{\zeta} + R^{N}) \leq C R^{1+\delta_0 + \zeta},
		\end{aligned}
	\end{equation*}
	where $C$ denotes a constant.
In view of H\"{o}lder's inequality, we obtain
	\begin{equation}
		\begin{aligned}
			\int_{B_R} |v - v_{x_0 , R}| \text{d}x & \leq C |B_R (x_0)|^{1 - \frac{1}{1+\delta_0}} R^{\frac{1+\delta_0 +\zeta}{1+\delta_0}} = C R^{N + \sigma},
		\end{aligned}
	\end{equation}
	which, combined with Campanato's Theorem \cite[Theorem 1.54]{maly}, provides the result. To prove the other items, we apply a similar argument, replacing the item $(a)$ of Theorem \ref{morreyg} by the respective item of the same result.
\end{proof}

\section{BMO regularity to the gradient and Log-Lipschitz estimates}
\hfill

This section is devoted to prove Theorem \ref{bmo1}.
\begin{proof}[Proof of Theorem \ref{bmo1}]
 Consider that $(\widetilde{a})$, $(\widetilde{b})$, or $(\widetilde{c})$ holds true.
 Let $B_r(x_0) \subset B_R(x_0)$ be open balls with $0 < r \leq R \leq R_0$, where $R_0$ is small. Consider $v^* \in W^{1,\Phi}(B_R(x_0))$ the $\Phi-$harmonic replacement of $v$ in $B_R(x_0),$ that is, $v^*$ satisfies \eqref{hreplace}.  Suppose also that $v^* = v$ in $\Omega \cap B_R(x_0)^{c}.$
Recall that 	 Lemma    \ref{comp} provides 
	\begin{equation}\label{comphar}
		\int_{B_r} \Phi(|\nabla v-(\nabla v)_{x_0 , r}|)\text{d}x\leq
			C\bigg(\frac{r}{R}\bigg)^{N+\sigma}\int_{B_R} \Phi(|\nabla v-(\nabla
			v)_{x_0 , R}|)\text{d}x +C\int_{B_R} \Phi(|\nabla v-\nabla v^*|)\text{d}x.
	\end{equation}

	The minimality of $v$, 	\eqref{m3.0},  the fact that $\| v^*\|_{L^{\infty}(B_R(x_0))} \leq \| v\|_{L^{\infty}(\Omega)},$ and Lemma \ref{f-properties} provides that
\begin{equation}\label{cs0}
\begin{aligned}
\int_{B_R} \Phi (|\nabla v|) - \Phi (|\nabla v^*|)\text{d}x& \leq CR^{N} +  \int_{\Omega_1 \cap B_R } g(v^{*} -v) \text{d}x+\int_{\Omega_2 \cap B_R } \overline{g}(v^{*} -v) \text{d}x+\int_{\Gamma \cap B_R } h (v-v^*)\text{d}x\\
&+ \int_{\Gamma \cap B_R } \rho'(Y((v^*)^{+}) - Y(v^{+})) \text{d}H^{N-1} + \int_{\Gamma \cap B_R } \widetilde{\rho}(Y((v^*)^{-}) - Y(v^{-})) \text{d}H^{N-1}.
\end{aligned}
\end{equation}
Since  $v^* -v \in W^{1,1}_{0}(B_R (x_0)),$ and $v^* = v$ in $\Omega \cap B_R(x_0)^{c}$ it follows that  $v^* -v \in W^{1,1}_{0}(\Omega)$. Thus, recalling \eqref{overb}, we have
\begin{equation}\label{cs1}
\begin{aligned}
	\int_{\Omega_1 \cap B_R } g (v^* -v)\text{d}x +  \int_{\Omega_2 \cap B_R } \overline{g} (v^* -v)\text{d}x  
	\leq  C\left( R^{N} + \varepsilon \int_{B_R} \Phi(|\nabla v^* - \nabla v|)\text{d}x \right),
\end{aligned}
\end{equation}
where $C>0$ is a constant, and $\varepsilon>0$ is arbitrary. 
From \eqref{bgood},\eqref{boundary1}-\eqref{boundary2},\eqref{cach1}-\eqref{c2} it follows that
\begin{equation}\label{cs3}
\int_{\Gamma \cap B_R }  \rho'(Y((v^*)^{+}) - Y(v^{+})) + \widetilde{\rho}(Y((v^*)^{-}) - Y(v^{-})) + h (v - v^*)\text{d}H^{N-1}\leq {C} \left(R^{N} + \varepsilon \int_{B_R} \Phi(|\nabla v - \nabla v^*|)\text{d}x\right),
\end{equation}
where ${C}>0$ is a constant, and $\varepsilon >0$ is arbitrary.

Consider that \eqref{ed} holds.
From item $(ii)$ of Lemma \ref{simon2} we have  $v \in W^{1,p}(\Omega),$ and $v^* \in W^{1,p}(B_R(x_0)).$
Thus, from  \eqref{bgood},\eqref{boundary1}-\eqref{boundary2},\eqref{cach1}-\eqref{c2}, and the Poincaré inequality  we obtain that
\begin{equation}\label{younged}
	\begin{aligned}
		\int_{\Gamma \cap B_R }  \rho'(Y((v^*)^{+}) - Y(v^{+})) &+\widetilde{\rho}(Y((v^*)^{-}) - Y(v^{-})) + h (v - v^*)\text{d}H^{N-1} \leq C \int_{\Gamma \cap B_R} |v-v^*|\text{d}H^{N-1}\\
		& \leq C\left( \int_{\Omega_1} |v -v^*| \text{d}x + \int_{\Omega_1} |\nabla v -\nabla v^*| \text{d}x  \right)\\
		& = C\left( \int_{\Omega_1 \cap B_R } |v -v^*| \text{d}x + \int_{\Omega_1 \cap B_R } |\nabla v -\nabla v^*| \text{d}x  \right) \\
		& \leq  C\left( \int_{ B_R} |v -v^*| \text{d}x + \int_{B_R} |\nabla v -\nabla v^*| \text{d}x  \right)\\
		& \leq C \int_{B_R} |\nabla v - \nabla v^*|\text{d}x,
	\end{aligned}
\end{equation}
where $C>0$ is a constant.
Since \eqref{ed} holds, it follows from \eqref{cs0}-\eqref{younged}, and Lemma  \ref{simon2} that

\begin{equation*}
	\begin{aligned}
		\int_{ B_R} |\nabla v -\nabla v^*|^{p}\text{d}x \leq C\int_{ B_R} \Phi (|\nabla v|) - \Phi(|\nabla v^*|)\text{d}x \leq  C\int_{ B_R} |\nabla v - \nabla v^*|\text{d}x,
	\end{aligned}
\end{equation*}
where $C >0$ is a constant.
From the previous estimate and the H\"{o}lder's inequality, we have 
\begin{equation*}
	\begin{aligned}
		\left( \int_{ B_R} |\nabla v - \nabla v^*|\text{d}x\right)^{p} \leq C R^{N(p-1)}\int_{ B_R} |\nabla v - \nabla v^*|^{p}\text{d}x 
		\leq C R^{N(p-1)}\int_{ B_R} |\nabla v - \nabla v^*|\text{d}x,
	\end{aligned}
\end{equation*}
which provides that
\begin{equation*}
	\int_{ B_R}|\nabla v - \nabla v^*| \text{d}x \leq CR^{N}.
\end{equation*}
The last inequality and Lemma \ref{comp} imply that
\begin{equation*}
	\begin{aligned}
		\int_{B_r} |\nabla v-(\nabla v)_{x_0,r}|\text{d}x & \leq
		C\bigg(\frac{r}{R}\bigg)^{N+\overline{\sigma}}\int_{B_R} |\nabla v-(\nabla
		v)_{x_0,R}|\text{d}x + C\int_{ B_R} |\nabla v-\nabla v^*|\text{d}x\\
		& \leq C\bigg(\frac{r}{R}\bigg)^{N+\overline{\sigma}}\int_{B_R} |\nabla v-(\nabla
		v)_{x_0,R}|\text{d}x  + CR^N,
	\end{aligned}
\end{equation*}
which combined with Lemma \ref{itg} imply
\begin{equation*}
	\int_{ B_r}|\nabla v-(\nabla v)_{x_0,r}|\text{d}x \leq Cr^N,
\end{equation*}
providing that  $\frac{\partial v}{\partial x_i} \in {\rm{BMO}}_{\rm loc}(\Omega ), i=1,\cdots,N$. Proceeding in the same way as in \cite[Section 5]{LQT} (see also \cite[Page 20]{piment} ), it follows the Log-Lipschitz regularity.

Suppose now that  $\Phi' = \phi$ is convex and that $\delta_0>1$. From $(i)$ of Lemma \ref{simon2} and \eqref{cs0}-\eqref{cs3}, we obtain
\begin{equation*}
\begin{aligned}
\int_{B_R} \Phi (|\nabla v -\nabla v^*|)\text{d}x \leq C \int_{B_R} \Phi (|\nabla v|) - \Phi(|\nabla v^*|)\text{d}x \leq C \left( R^{N} + \varepsilon \int_{B_R} \Phi (|\nabla v - \nabla v^*|)\right)\text{d}x,
\end{aligned}
\end{equation*}
where $C>0$ is a constant, and $\varepsilon>0$ is arbitrary.
A suitable choice of $\varepsilon$ in the last inequality implies that
\begin{equation*}
\int_{B_R}\Phi(|\nabla v -\nabla v^*|)\text{d}x \leq C R^{N},
\end{equation*}
which combined with  \eqref{comphar}, provides that
\begin{equation}\label{iterar}
\int_{B_r} \Phi(|\nabla v - (\nabla v)_{x_0, r}|)\text{d}x \leq C\left( \dfrac{r}{R}\right)^{N + \sigma} \int_{B_R} \Phi(|\nabla v - (\nabla v)_{x_0, R}|)\text{d}x +CR^{N},
\end{equation}
for all $0<r \leq R,$ where $C>0$ is a constant. 
As pointed out in \cite[Page 36]{zhengt}, the function $${\psi(r):= \displaystyle\int_{B_r} \Phi(|\nabla v - (\nabla v)_{x_0, r}|)\text{d}x }$$ satisfies the hypotheses of Lemma \ref{itg}. Thus, from \eqref{iterar}, and Lemma \ref{itg} we obtain that

\begin{equation}\label{fim1}
\int_{B_r} \Phi (|\nabla v - (\nabla v)_{x_0,r}|) \text{d}x\leq C r^N,
\end{equation}
for all $0<r\leq R\leq R_0,$ where $C>0$ is a constant that does not depend on $r.$
We claim that there is a constant $K>0$ such that 
\begin{equation*}
\fint_{B_r } |\nabla v - (\nabla v)_{x_0 , r}|\text{d}x \leq K,
\end{equation*}
for all $0<r\leq R, $ where $K>0$ is a constant that does not depend on $r,R,$ and $x_0 .$ 
In fact, if $$r^{-N} \displaystyle\int_{B_r} \Phi (|\nabla v - (\nabla v)_{x_0 , r}|)\text{d}x \leq 1, \mbox{ then } \displaystyle\int_{B_r} \Phi (|\nabla v - (\nabla v)_{x_0 , r}|)\text{d}x \leq r^N.$$
Suppose that $$r^{-N} \displaystyle\int_{B_r} \Phi (|\nabla v - (\nabla v)_{x_0 , r}|)\text{d}x > 1.$$
 Computing the derivative of the function $\Psi(t) := \Phi(t) - \Phi(1)t$, for $ t\geq1, $ and using Lemma \ref{G-properties}, items $(\phi_1)$ and $(\Phi_2)$, we conclude that $\Psi$ is nondecreasing with $\Psi(1)=0$. Then  $\Psi(t) \geq 1,$ for all $t\geq 1$, which provides that $\Phi(t) \geq \Phi(1)t,$ for all $t\geq 1.$
Therefore,
\begin{equation}\label{quase0}
\Phi\left( r^{-N} \int_{B_r} |\nabla v - (\nabla v)_{x_0 ,r}|\text{d}x\right)\geq \Phi(1) \left( r^{-N} \int_{B_r} |\nabla v - (\nabla v)_{x_0 ,r}| \text{d}x\right).
\end{equation}
From \eqref{fim1}-\eqref{quase0}, and the convexity of $\Phi,$ which is provided by item $(\Phi_1)$ of Lemma \ref{G-properties}, it follows that
\begin{equation*}
\begin{aligned}
\Phi(1) \left( r^{-N} \int_{B_r} |\nabla v - (\nabla v)_{x_0 ,r}| \text{d}x\right) & \leq \Phi \left( r^{-N} \int_{B_r} |\nabla v - (\nabla v)_{r}| \text{d}x\right)\leq \frac{C}{r^N} \int_{B_r} \Phi (|\nabla v - (\nabla v)_{x_0 , r}|)\text{d}x \leq C,
\end{aligned}
\end{equation*}
for some constant $C>0$ that does not depend on $r$, which proves the claim. Arguing as before we obtain the Log-Lipschitz regularity. 
\end{proof}

\section{The case $h \in W^{1,p^{\prime} + \varepsilon}(\Omega) $ }

Below it will be presented the proof of Theorem \ref{last1}.
\begin{proof}[Proof of Theorem \ref{last1}]
Consider an open ball $B_R(x_0) \subset \Omega$ and let $v^* \in W^{1,\Phi}(B_R(x_0))$ be  the $\Phi-$harmonic replacement of the minimizer $v$ in $B_R(x_0),$ that is, $v^*$ satisfies \eqref{hreplace}. Suppose also that $v^* = v$ in $\Omega \cap B_R(x_0)^{c}.$ Thus, the estimate \eqref{m3.0} holds. The idea is to estimate each term on the right-hand side of \eqref{m3.0}. From the continuous  embedding $W^{1,1}(\Omega_1)\hookrightarrow L^{1}(\Gamma)$ it follows that
\begin{equation*}
\begin{aligned}
\int_{\Gamma \cap B_R} |h||v-v^*|\text{d}H^{N-1} &\leq \int_{\Gamma}|h||v-v^*|\text{d}H^{N-1}\leq C \left( \int_{B_R} |h||v-v^*|\text{d}x + \int_{B_R} |\nabla (h (v-v^*))|\text{d}x \right)\\
&\leq C  \left( \int_{B_R}  |h||v-v^*|\text{d}x +\int_{B_R}  |\nabla h||v-v^*| \text{d}x + \int_{B_R} |h| |\nabla v - \nabla v^*| \text{d}x\right)\\
& := C(I_1 + I_2 + I_3).
\end{aligned}
\end{equation*}	
Since $v-v^* \in L^{1+\delta_0}(B_R)$, by Poincaré's Inequality,  we obtain that

\begin{equation}\label{i1}
\begin{aligned}
I_1  &\leq  
\left( \int_{B_R} |h|^{p^{\prime} + \varepsilon} \text{d}x \right)^{\frac{1}{p^{\prime} + \varepsilon}}  
\left( \int_{B_R} |v-v^*|^{\frac{p^{\prime} + \varepsilon}{p^{\prime} + \varepsilon - 1}} \text{d}x \right)^{\frac{p^{\prime} + \varepsilon - 1}{p^{\prime} + \varepsilon}}\\
&\leq C  
\left( \left( \int_{B_R} |v-v^*|^{\frac{p^{\prime} + \varepsilon}{p^{\prime} + \varepsilon - 1} \cdot \frac{p^{\prime} + \varepsilon - 1}{p^{\prime} + \varepsilon} p} dx \right)^{\frac{p^{\prime} + \varepsilon}{p^{\prime} + \varepsilon - 1} \cdot \frac{1}{p}}  
\left(\int_{B_R} \text{d}x\right)^{\frac{\frac{p^{\prime} +\varepsilon -1}{p^{\prime} +\varepsilon} p-1}{\frac{p^{\prime} +\varepsilon -1}{p^{\prime} +\varepsilon}p}} \right)^{\frac{p^{\prime}+\varepsilon -1}{p^{\prime}+ \varepsilon} }\\
&\leq C R^{N \frac{p^{\prime} + \varepsilon - 1}{p^{\prime} + \varepsilon} \frac{p-1}{p}} \|v-v^*\|_{L^{p}(B_R)} \leq C R^{N \frac{p^{\prime} + \varepsilon - 1}{p^{\prime} + \varepsilon} \frac{p-1}{p}} \|\nabla (v-v^*)\|_{L^{p}(B_R)}.
\end{aligned}
\end{equation}

Note also that 

\begin{equation}\label{i2}
	\begin{aligned}
		I_2  &\leq  
		\left( \int_{B_R} |\nabla h|^{p^{\prime} + \varepsilon} dx \right)^{\frac{1}{p^{\prime} + \varepsilon}}  
		\left( \int_{B_R} |v-v^*|^{\frac{p^{\prime} + \varepsilon}{p^{\prime} + \varepsilon - 1}} dx \right)^{\frac{p^{\prime} + \varepsilon - 1}{p^{\prime} + \varepsilon}}\\
		&\leq C  
		\left( \left( \int_{B_R} |v-v^*|^{\frac{p^{\prime} + \varepsilon}{p^{\prime} + \varepsilon - 1} \cdot \frac{p^{\prime} + \varepsilon - 1}{p^{\prime} + \varepsilon} p} dx \right)^{\frac{p^{\prime} + \varepsilon}{p^{\prime} + \varepsilon - 1} \cdot \frac{1}{p}}  
		\left(\int_{B_R} \text{d}x\right)^{\frac{\frac{p^{\prime} +\varepsilon -1}{p^{\prime} +\varepsilon} p-1}{\frac{p^{\prime} +\varepsilon -1}{p^{\prime} +\varepsilon}p}} \right)^{\frac{p^{\prime}+\varepsilon -1}{p^{\prime}+ \varepsilon} }\\
		&\leq C R^{N \frac{p^{\prime} + \varepsilon - 1}{p^{\prime} + \varepsilon} \frac{p-1}{p}} \|v-v^*\|_{L^{p}(B_R)} \leq C R^{N \frac{p^{\prime} + \varepsilon - 1}{p^{\prime} + \varepsilon} \frac{p-1}{p}} \|\nabla (v-v^*)\|_{L^{p}(B_R)},
	\end{aligned}
\end{equation}
and 

\begin{equation}\label{i3}
	\begin{aligned}
		I_3  &\leq  
		\left( \int_{B_R} | h|^{p^{\prime} + \varepsilon} \text{d}x \right)^{\frac{1}{p^{\prime} + \varepsilon}}  
		\left( \int_{B_R} |\nabla v-\nabla v^*|^{\frac{p^{\prime} + \varepsilon}{p^{\prime} + \varepsilon - 1}} \text{d}x \right)^{\frac{p^{\prime} + \varepsilon - 1}{p^{\prime} + \varepsilon}}\\
		&\leq C  
		\left( \left( \int_{B_R} |\nabla v-\nabla v^*|^{\frac{p^{\prime} + \varepsilon}{p^{\prime} + \varepsilon - 1} \cdot \frac{p^{\prime} + \varepsilon - 1}{p^{\prime} + \varepsilon} p} \text{d}x \right)^{\frac{p^{\prime} + \varepsilon}{p^{\prime} + \varepsilon - 1} \cdot \frac{1}{p}}  
		\left(\int_{B_R} \text{d}x\right)^{\frac{\frac{p^{\prime} +\varepsilon -1}{p^{\prime} +\varepsilon} p-1}{\frac{p^{\prime} +\varepsilon -1}{p^{\prime} +\varepsilon}p}} \right)^{\frac{p^{\prime}+\varepsilon -1}{p^{\prime}+ \varepsilon} }\\
		& \leq C R^{N \frac{p^{\prime} + \varepsilon - 1}{p^{\prime} + \varepsilon} \frac{p-1}{p}} \|\nabla (v-v^*)\|_{L^{p}(B_R)}.
	\end{aligned}
\end{equation}
Since $g\in L^{q}(\Omega_1), \overline{g}\in L^{\overline{q}}(\Omega_2),$ with  $\min\{q,\overline{q}\} \geq N >p^{\prime} + \varepsilon,$ we have   $g \in L^{p^{\prime} + \varepsilon}(\Omega_1),$ and $\overline{g} \in L^{p^{\prime} + \varepsilon}(\Omega_2)$. Thus, by applying the reasoning of \eqref{i1}-\eqref{i3} we obtain that
\begin{equation}\label{g1}
\int_{\Omega_1 \cap B_R} |g||v^*-v| \text{d}x \leq C R^{N \frac{p^{\prime} + \varepsilon - 1}{p^{\prime} + \varepsilon} \frac{p-1}{p}} \|\nabla (v-v^*)\|_{L^{p}(B_R)},
\end{equation}
and 
\begin{equation}\label{g2}
	\int_{\Omega_2 \cap B_R} |\overline{g}||v^*-v| \text{d}x \leq C R^{N \frac{p^{\prime} + \varepsilon - 1}{p^{\prime} + \varepsilon} \frac{p-1}{p}} \|\nabla (v-v^*)\|_{L^{p}(B_R)}.
\end{equation}
Since $y=Y^{\prime}$ is bounded, we obtain from \eqref{boundary1}-\eqref{boundary2},  the embedding $W^{1,1}(\Omega_1)\hookrightarrow L^{1}(\Gamma),$ and the reasoning of \eqref{i1}-\eqref{i3}  that
\begin{equation}\label{ybounded}
\begin{aligned}
	\int_{\Gamma \cap B_R} \rho'(Y((v^*)^{+}) - Y(v^{+})) \text{d}H^{N-1} &+ 	\int_{\Gamma \cap B_R} \widetilde{\rho}(Y((v^*)^{-}) - Y(v^{-})) \text{d}H^{N-1} \leq C \int_{\Gamma \cap B_R} | v^* - v| \text{d}H^{N-1}\\
	&\leq C \int_{\Gamma } | v^* - v| \text{d}H^{N-1}\leq C \left( \int_{\Omega_1} |v -v^*| \text{d}x + \int_{\Omega_1} |\nabla v - \nabla v^*|\right)\\
	&\leq C \left( \int_{B_R} |v -v^*| \text{d}x + \int_{B_R} |\nabla v - \nabla v^*|\text{d}x\right)\\
	& \leq C R^{N \frac{p^{\prime} + \varepsilon - 1}{p^{\prime} + \varepsilon} \frac{p-1}{p}} \|\nabla (v-v^*)\|_{L^{p}(B_R)}.
\end{aligned}
\end{equation}
The boundedness of  $f=F^{\prime}$ combined with the Mean Value Theorem, and Lemma \ref{ineqpart}, implies  that
\begin{equation*}
\begin{aligned}
\int_{\Omega_1 \cap B_R} \theta'(F((v^*)^{+}) - F(v^{+})) \text{d}x &+  \int_{\Omega_1 \cap B_R }  \widetilde{\theta}(F((v^*)^{-}) - F(v^{-})) \text{d}x \leq C \int_{B_R} |(v^*)^{+} -v^{+}|\text{d}x \\
& \leq C \int_{B_R} |v^* -v|\text{d}x \leq C  R^{N \frac{p^{\prime} + \varepsilon - 1}{p^{\prime} + \varepsilon} \frac{p-1}{p}} \|\nabla (v-v^*)\|_{L^{p}(B_R)},\\
\end{aligned}
\end{equation*}
for some constant $C>0$, which implies

\begin{equation}\label{casew1}
\begin{aligned}
\int_{\Omega_1 \cap B_R
} W(v^*) - W(v)\text{d}x & \leq C R^{N \frac{p^{\prime} + \varepsilon - 1}{p^{\prime} + \varepsilon} \frac{p-1}{p}} \|\nabla (v-v^*)\|_{L^{p}(B_R)},\\
\end{aligned}
\end{equation}
where $C>0$ is a   constant. Similarly, we obtain that
\begin{equation}\label{w2}
\begin{aligned}
	\int_{\Omega_1 \cap B_R} \overline{W}(v^*) - \overline{W}(v)\text{d}x & \leq C R^{N \frac{p^{\prime} + \varepsilon - 1}{p^{\prime} + \varepsilon} \frac{p-1}{p}} \|\nabla (v-v^*)\|_{L^{p}(B_R)},
\end{aligned}
\end{equation}
where $C>0$ is a   constant.
Thus,  from \eqref{m3.0}, \eqref{i1}-\eqref{w2} it follows that
\begin{equation}\label{quase}
\begin{aligned}
\int_{B_R} \Phi(|\nabla v|) \text{d}x - \int_{B_R} \Phi (|\nabla v^*|) \text{d}x &\leq C R^{N \frac{p^{\prime} + \varepsilon - 1}{p^{\prime} + \varepsilon} \frac{p-1}{p}} \|\nabla (v-v^*)\|_{L^{p}(B_R)},
\end{aligned}
\end{equation}
where $C >0$ is a  constant.
Since  \eqref{ed} holds, we have  from Lemma \ref{simon2},  and \eqref{quase} that

\begin{equation}\label{quasee}
\begin{aligned}
\int_{B_R} |\nabla v - \nabla v^*|^{p}  \text{d}x &\leq C  \int_{B_R} \Phi(|\nabla v|) \text{d}x - \int_{B_R} \Phi (|\nabla v^*|) \text{d}x \leq C R^{N \frac{p^{\prime} + \varepsilon - 1}{p^{\prime} + \varepsilon} \frac{p-1}{p}} \|\nabla (v-v^*)\|_{L^{p}(B_R)}.
\end{aligned}
\end{equation}
Note that \eqref{quasee}, combined with
$
\frac{p^{\prime} + \varepsilon - 1}{p^{\prime} + \varepsilon} \frac{p-1}{p} = \frac{\varepsilon (p-1)}{p(p^{\prime} + \varepsilon)},
$
implies that
\begin{equation}\label{fimm}
\int_{B_R} |\nabla v - \nabla v^*|^{p} \leq C R^{N\frac{\varepsilon}{p^{\prime}+\varepsilon}},
\end{equation}
for some constant $C>0.$
We claim that 
\begin{equation}\label{claimalmost}
\int_{B_R} |\nabla v - \nabla v^*|\text{d}x \leq C R^{N \frac{\varepsilon}{p^{\prime} + \varepsilon}},
\end{equation}
for some constant $C>0.$ In fact, if $R^{- N \frac{\varepsilon}{p^{\prime} + \varepsilon}} \cdot \displaystyle\int_{B_R} |\nabla v - \nabla v^*|\text{d}x \leq 1,$ then there is nothing to prove. 
Suppose that $ R^{-N \frac{\varepsilon}{p^{\prime} + \varepsilon}} \cdot\displaystyle\int_{B_R} |\nabla v - \nabla v^*|\text{d}x >1.$ From the fact that $t \leq t^{p},$ for all $t \geq 1,$ and the Jensen's inequality, and \eqref{fimm} we obtain 
\begin{equation*}
\begin{aligned}
R^{- N \frac{\varepsilon}{p^{\prime} + \varepsilon}}\int_{B_R} |\nabla v - \nabla v^*|\text{d}x  & \leq \left( R^{- \frac{N\varepsilon}{p^{\prime}+\varepsilon}}\int_{B_R} |\nabla v - \nabla v^*| \text{d}x \right)^{p} = R^{-N \frac{\varepsilon p}{p^{\prime} + \varepsilon}} |B_R |^{p} \left( \fint_{B_R} |\nabla v - \nabla v^*| \text{d}x\right)^{p}\\
& \leq R^{-N \frac{\varepsilon p}{p^{\prime} + \varepsilon}} |B_R|^{p} \frac{1}{|B_R|}\int_{B_R} |\nabla v - \nabla v^*|^{p} \text{d}x \\
&=CR^{-N \frac{\varepsilon p}{p^{\prime} + \varepsilon}} R^{N(p-1)}\int_{B_R} |\nabla v - \nabla v^*|^{p} \text{d}x \\
& \leq C R^{- N \frac{\varepsilon p}{p^{\prime} + \varepsilon}} R^{N(p-1)}R^{\frac{N\varepsilon}{p^{\prime}+\varepsilon}} =C R^{\frac{Np}{p^{\prime} + \varepsilon}}  \leq  {C (\text{diam}(\Omega))^{\frac{N p}{p^{\prime} + \varepsilon}}}.
\end{aligned}
\end{equation*} 
Thus,  the claim \eqref{claimalmost} is proved.
From Lemma \ref{comp}, and \eqref{claimalmost} it follows that

\begin{equation*}
\begin{aligned}
\int_{B_r} |\nabla v-(\nabla v)_{x_0,r}|\text{d}x& \leq
C\bigg(\frac{r}{R}\bigg)^{N+\sigma}\int_{B_R} |\nabla v-(\nabla
v)_{x_0,R}|\text{d}x +C\int_{B_R} |\nabla v-\nabla v^*|\text{d}x\\
& \leq C\bigg(\frac{r}{R}\bigg)^{N+\sigma}\int_{B_R} |\nabla v-(\nabla
v)_{x_0,R}|\text{d}x +CR^{N\frac{\varepsilon}{p^{\prime} + \varepsilon}},
\end{aligned}
\end{equation*}
which, combined  with Lemma \ref{itg}, provides
\begin{equation*}
\int_{B_r} |\nabla v-(\nabla v)_{x_0,r}|\text{d}x\leq Cr^{N \frac{\varepsilon}{p^{\prime} + \varepsilon}},
\end{equation*}
for all $0 < r \leq R. $ The last inequality  and the hypothesis $p^{\prime} \leq \varepsilon $  provide that
\begin{equation*}
\begin{aligned}
r^{- \left( N - \frac{N \varepsilon}{p^{\prime} + \varepsilon}\right)} \int_{B_r} |\nabla v-(\nabla v)_{x_0,r}|\text{d}x & \leq C r^{\frac{2N\varepsilon}{p^{\prime}+\varepsilon} - N}  \leq C,
\end{aligned}
\end{equation*}
for some constant $C>0.$
 Therefore $\frac{\partial v}{\partial x_i} \in L^{1, N - N \frac{\varepsilon}{p^{\prime} + \varepsilon}}_{C}(\Omega^{\prime}),$ for all domain $\Omega^{\prime} \subset \subset \Omega,$ and $i=1,\cdots,N.$ 
 Regarding the remaining part of the result, note that if $\phi(t) = t^{p-1} ,$ we obtain from Lemma \ref{comp}, and  \eqref{fimm} that
 \begin{equation*}
 \begin{aligned}
 \int_{B_r} |\nabla v-(\nabla v)_{x_0,r}|^{p}\text{d}x & \leq
 C\bigg(\frac{r}{R}\bigg)^{N+ \sigma}\int_{B_R} |\nabla v-(\nabla
 v)_{x_0,R}|^p\text{d}x+C\int_{B_R} |\nabla v-\nabla v^*|^p\text{d}x\\
 & \leq C\bigg(\frac{r}{R}\bigg)^{N+\sigma}\int_{B_R}|\nabla v-(\nabla
 v)_{x_0,R}|^p\text{d}x + CR^{N\frac{\varepsilon}{p^{\prime} + \varepsilon}},
\end{aligned}
 \end{equation*}
 which combined with Lemma \ref{itg} provides
 
 \begin{equation}\label{campfinal}
 \int_{B_r} |\nabla v-(\nabla v)_{x_0,r}|^{p}\text{d}x  \leq C r^{N\frac{\varepsilon}{p^{\prime} + \varepsilon}},
 \end{equation} 
 for all $0 <r \leq R.$ Since $p^{\prime}\leq \varepsilon$ we have from \eqref{campfinal} that
 
 \begin{equation*}
 	\begin{aligned}
 		r^{- \left( N - \frac{N \varepsilon}{p^{\prime} + \varepsilon}\right)} \int_{B_r} |\nabla v-(\nabla v)_{x_0,r}|^{p}\text{d}x & \leq C r^{\frac{2N\varepsilon}{p^{\prime}+\varepsilon} - N}  \leq C,
 	\end{aligned}
 \end{equation*}
which implies that $\tfrac{\partial v}{\partial x_i} \in L^{p, \lambda}_{C , \rm{loc}}(\Omega),$  $i=1,\cdots,N,$ where $\lambda:= N - N \frac{\varepsilon}{p^{\prime} + \varepsilon} = N(1 - \frac{\varepsilon}{p^{\prime} + \varepsilon})= \frac{Np}{p +p\varepsilon - \varepsilon}.$ 
 The residual argument follows from the embeddings of the Campanato's spaces $L^{p, \lambda}_{C , \rm{loc}}(\Omega)$ on the H\"older space $C_{\rm loc}^{0,\alpha}(\Omega)$, for the details, see the proof of \cite[Theorem 3]{piment}.
 \end{proof}

\section*{Declarations}

\subsection*{Authors contributions}
All authors of this manuscript contributed equally to this work.
\subsection*{Data Availability}
No data were used to support this study.
\subsection*{Conflicts of Interest}
The authors declare that there are no conflicts of interest.
\subsection*{Use of AI tools declaration}
The authors declare they have not used Artificial Intelligence (AI) tools in the creation of this article.
\subsection*{Funding}
 R. Lima was partially supported by CNPq/Brazil 306.411/2022-9 and 176.596/2023-2

\vspace{1cm}

\noindent\textbf{Romildo Lima}\\
Unidade Acadêmica de Matemática\\
Universidade Federal de Campina Grande \\
58429-970, Campina Grande - PB, Brazil\\
romildo.nascimento@professor.ufcg.edu.br

\vspace{0.5cm}

\noindent\textbf{Mayra Soares}\\
Departamento de Matemática\\
	Universidade de Brasília \\
	70910-900, Brasília - DF, Brazil\\
mayra.soares@unb.br

\vspace{0.5cm}

\noindent\textbf{Leandro Tavares}\\
Centro de Matemática, Computação e Cognição\\
Universidade Federal do ABC, \\
63048-080, Santo André - SP, Brazil\\
leandrolstav@gmail.com

\end{document}